\documentclass[preprint]{elsarticle}
\usepackage{times,epsf,amsmath,amssymb,cite}
\usepackage{caption2}
\usepackage{hyperref}
\usepackage{amsthm}
\usepackage{graphicx} 
\usepackage{mathabx}
\usepackage{makecell}
\usepackage{wrapfig}

\newcommand{\RNum}[1]{\uppercase\expandafter{\romannumeral #1\relax}}

\newtheorem{Theorem}{Theorem}
\newtheorem*{theorem*}{Theorem}
\newtheorem{Proposition}{Proposition}
\newtheorem{Corollary}{Corollary}
\newtheorem{Definition}{Definition}

\usepackage{bbm}
\usepackage{bm}
\newcommand\norm[1]{\left\lVert#1\right\rVert}
\newcommand\snorm[1]{\left\vert#1\right\vert}
\newcommand{\e}{\mathrm{e}}
\newcommand\bb[1]{\mathbf{#1}}

\newcommand{\R}{\mathbb{R}}
\newcommand{\N}{\mathbb{N}}
\newcommand{\C}{\mathcal{C}}
\newcommand{\M}{\mathcal{M}}
\newcommand{\Yy}{\mathcal{Y}}
\newcommand{\Ls}{\mathrm{L}}

\newcommand{\D}{\mathcal{D}}
\newcommand{\Dd}{\mathrm{D}}
\newcommand{\K}{\mathrm{K}}
\newcommand{\DD}{\mathrm{D}\otimes\mathrm{D}}
\newcommand{\DDD}{\mathrm{D}^2\otimes\mathrm{D}^2}

\newcommand{\I}{\mathrm{I}}
\newcommand{\LL}{\Ls_1\otimes\Ls_2}
\newcommand{\B}{\mathcal{B}}
\newcommand{\w}{\mathrm{weak}^{\star}}

\journal{Applied and Computational Harmonic Analysis}
\begin{document}
\begin{frontmatter}

\title{Variational Tensor-Product Splines}

 \author[1]{Vincent Guillemet\corref{cor1}}
 \ead{vincent.guillemet@epfl.ch}
 
 \author[1]{Michaël Unser}
 \ead{michael.unser@epfl.ch}

 \cortext[cor1]{Corresponding author}

 \affiliation[1]{organization={Biomedical Imaging Group,
École polytechnique fédérale de Lausanne (EPFL)},
 addressline={Station 17},
 city={Lausanne},
 postcode={CH-1015},
 state={Vaud},
 country={Suisse}}

\begin{abstract}
Multidimensional continuous-domain inverse problems are often solved by the minimization of a loss functional, formed as the sum of a data fidelity and a regularization.
In this work, we present a new construction where the regularization is itself built as the sum of two terms: i) the $\M$ norm of the regularizing operator $\Ls_1\otimes\Ls_2$, with $\Ls_1$ and $\Ls_2$ being two one-dimensional differential operators; ii) a bounded-variation norm that regularizes on the infinite-dimensional nullspace of $\Ls_1\otimes\Ls_2$. In this construction, we show that the extreme points of the solution set are the tensor product of one-dimensional splines, with a number of atoms upper-bounded in term of the number of data points. Further, when the data of the inverse problem is localized, we reveal that the term ii) must take the form of a sum of bounded-variation norms, precomposed with partial derivative of different orders.

\end{abstract}


\end{frontmatter}

%%%%%%%%%%%%%%%%%%%%%%%%%%%%%%%%%%%%%%%%%%%%%%%%%%%%%
\section{Introduction} 
Variational splines, bound to an ordinary differential operator (ODO) $\Ls$, are fundamental tools in approximation theory and signal processing \citep{de1972calculating,schumaker2007spline,unser1999splines}. The key properties of these functions are: a) support optimality \citep{delgado2012exponential}; b) optimal approximation accuracy \citep{Moms}; c) hierarchical multiresolution \citep{burt1983multiresolution,mallat1989theory}; d) simplicity of implementation \citep{de1978practical}; e) optimal rates for nonlinear approximation \citep{mallat1999wavelet}. Although used initially for interpolation and smoothing problems \citep{reinsch1967smoothing,schoenberg1973cardinal}, these splines have led to the successful development of algorithms for the resolution of general, one-dimensional, continuous-domain inverse problems (IPs) \citep{debarre2022uniqueness,debarre2019b,fisher1975spline,unser2017splines}. In parallel, the surge of imaging techniques, from biomedical \citep{bertero2021introduction} to radio-astronomy \citep{thompson2017interferometry} applications, lead to an extension of such signal-processing methods to multiple dimensions \citep{lim1990two}. In the variational formulations, the unknown image is recovered as the minimum of a loss functional, itself formed as the sum of a data-fidelity term, and of a regularization term that promotes specific types of spline solutions.

\subsection{Multidimensional Variational Splines}
Our interest lies in the sparsity-promoting variational formulations that make use of the $\M$ norm \citep[Chapter 6]{rudin} to promote solutions that are linear combinations of a small number of atoms taken in a dictionary \citep{bredies2020sparsity,flinth2019exact,unser2017splines}. When these atoms are spline-basis functions, whose properties are determined by a unique partial differential operator (PDO) $\Ls$ \citep[Theorem IX.23]{reed2003methods}, the choice of $\Ls$ becomes critical. Relevant examples include radial-basis functions through polyharmonic splines, associated to the fractional Laplacian $\Ls=(-\Delta)^{\gamma/2}$ \citep{duchon1977splines, madych1990polyharmonic} or to the Sobolev operator $\Ls=(1+\Delta)^{\gamma/2}$ \citep{ward2014approximation}. Nevertheless, the support of those splines is larger than needs be and their implementation is difficult.

As the ambient dimension grows higher, the finding of splines and associated PDO with Properties a)-e) becomes challenging: i) it is in general impossible to construct B-spline-like basis functions that are compactly supported; ii) it is necessary to use more complicated PDO, of higher order, for the Green's function to be continuous or  integrable; iii) the null space of $\Ls$ is often infinite-dimensional, which then requires a specialized treatment.

One class of PDO that mitigate Issues i)-iii) is the tensor product of ODO, with 
\begin{equation}
    \Ls=\LL=(\Dd-\alpha_1\mathrm{I})^{N_1}\otimes(\Dd-\alpha_2\mathrm{I})^{N_2},
\end{equation}
where the operators $\Ls_1,\Ls_2$ encode the underlying expected regularity, in the $x_1,x_2$ dimensions. Although one could argue that this technique yields the canonic multi-dimensional extension of signal processing, surprisingly, tensor splines have no variational interpretation. The primary difficulty there is that the null space of the natural regularization $\norm{[\LL]\{\cdot\}}_{\M(\R^2)}$ is infinite-dimensional.

\subsection{Contributions: Tensor Product Variational Splines}
\label{sec:1.4}
Our contribution lies in the definition of $[\LL]$-splines and the presentation of a new variational framework that promotes such functions. This framework revolves around the optimization problem in
\begin{equation}
\label{eq:3.1.30}
\mathcal{V}=\underset{f\in\mathcal{X}}{\text{argmin}}\quad\mathcal{J}(f)=\underset{f\in\mathcal{X}}{\text{argmin}}\quad \left(E(\bb{y},\langle f,\bm{\nu}\rangle)+\lambda\snorm{f}_{\M_{\LL}}\right),
\end{equation}
with the novel seminorm $\snorm{\cdot}_{\M_{\LL}}$ (see \eqref{def:seminorm}) that induces the following properties.
\begin{itemize}
    \item The native space $\mathcal{X}$ of $[\LL]$-splines is either $\M_{\LL}(\R^2)$, or $\M_{\LL}(\bb{K})$, where $\bb{K}=\mathrm{K}_1\times\mathrm{K}_2$ is a rectangle. Aside from considerations on the null space, this native space is the space of distributions whose $[\LL]$-weak derivative is a two-dimensional bounded measure supported in $\R^2$ ($\bb{K}$, respectively). 
    \item The seminorm $\snorm{\cdot}_{\M_{\LL}}:\M_{\LL}\to\R^+$, defined in \eqref{def:seminorm}, includes specialized regularization terms that make its nullspace finite-dimensional and equal to $\mathcal{N}_{\Ls_1}\otimes\mathcal{N}_{\Ls_2}$.
    \item The operator $\langle\cdot,\bm{\nu}\rangle:\M_{\Ls_1\otimes\Ls_2}(\R^2)\to\R^{M}$ is the forward measurement operator. For the optimization problem to be well-posed with a guaranteed existence of solutions, the component $\nu_m$ of $\bm{\nu}=(\nu_m)_{m=1}^M$ must belong to some appropriate vector space, noted $\C_{\LL}(\R^2)$, which limits their degree of singularity.
\end{itemize}
The vector $\bb{y}\in\R^M$ represents the available measurements. The functional $E$ is the data-fidelity functional that measures the discrepancy between the data $\bb{y}$ and the simulated measurement $\langle f,\bm{\nu}\rangle$. The technical requirements are listed in Assumptions 1 to 6, given in Section \ref{sec:5.1}. 

\subsection{Contributions: New Representer Theorem}
For $i\in\{1,2\}$, we denote by  $g_{\Ls_i}$ the causal Green's function of $\Ls_i$ and by $\mathcal{N}_{\Ls_i}=\text{span}\{p_{i,n}\}_{n=1}^{N_i}$ its null space, of dimension $N_i$. Our main result, Theorem \ref{th:RTcompact}, can be summarized as follows.

\begin{figure}
    \centering
    \includegraphics[width=0.55\linewidth]{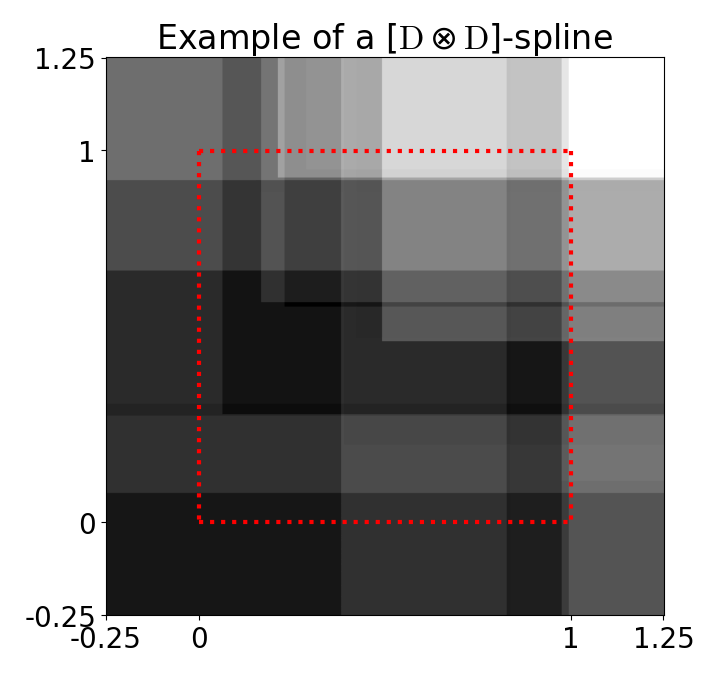}
    \caption{Example of a $[\Dd\otimes\Dd]$-spline   of the form \eqref{eq:1.7} with $K_1=5,K_2=5,$
    and $K=10$. The dotted square 
    represents the boundary 
    of the unit square $[0,1]^2$. The knots and the amplitudes have been chosen randomly. The plot is made on the extended square $[-0.25,1.25]^2$ to display the canonic behaviour of the spline outside of its supporting domain $[0,1]^2.$ On the $y$ axis ($x$ axis, respectively) one can observe the localisation of the knots $z_{m'}$ (the knots $y_m$, respectively) from the colour gradient. Finally the corners inside the dotted square are the localizations of the knots $(x_{1,k},x_{2,k})$.}
    \label{fig:1}
\end{figure}

\emph{For $\mathcal{X}=\M_{\LL}(\bb{K})$, the solution set $\mathcal{V}$ is nonempty, compact, and convex. Its extreme points are $[\LL]$-splines of the form
\begin{align}
\label{eq:1.7}
f(t_1,t_2)=&\sum_{k=1}^Ka_{k}g_{\Ls_1}(t_1-x_{1,k}) g_{\Ls_2}(t_2-x_{2,k})+\sum_{n=1}^{N_1}\sum_{m=1}^{K_{1,n}}b_{n,m}p_{1,n}(t_1) g_{\Ls_2}(t_2-y_{n,m})\nonumber\\
&+\sum_{n'=1}^{N_2}\sum_{m'=1}^{K_{2,n'}}c_{n',m'}g_{\Ls_1}(t_1-z_{n',m'}) p_{2,n'}(t_2)+\sum_{n,n'=1,1}^{N_1,N_2}d_{n,n'}p_{1,n}(t_1) p_{2,n'}(t_2).
\end{align}    
The parameters of \eqref{eq:1.7} satisfy 
\begin{itemize}
    \item []\emph{Sparsity:} The number of spline atoms $\left(K+\sum_{n=1}^{N_1}K_{1,n}+\sum_{n=1}^{N_2}K_{2,n}\right)$ is upper-bounded by $(M-N_1N_2)$.
    \item []\emph{Localization:} The knots are localized as $(x_{1,k},x_{2,k})\in\bb{K},y_{n,m}\in\mathrm{K}_2$, and $z_{n',m'}\in\mathrm{K}_1.$ 
\end{itemize}}

\begin{figure}[h!]
    \centering
    \includegraphics[width=1.2\linewidth]{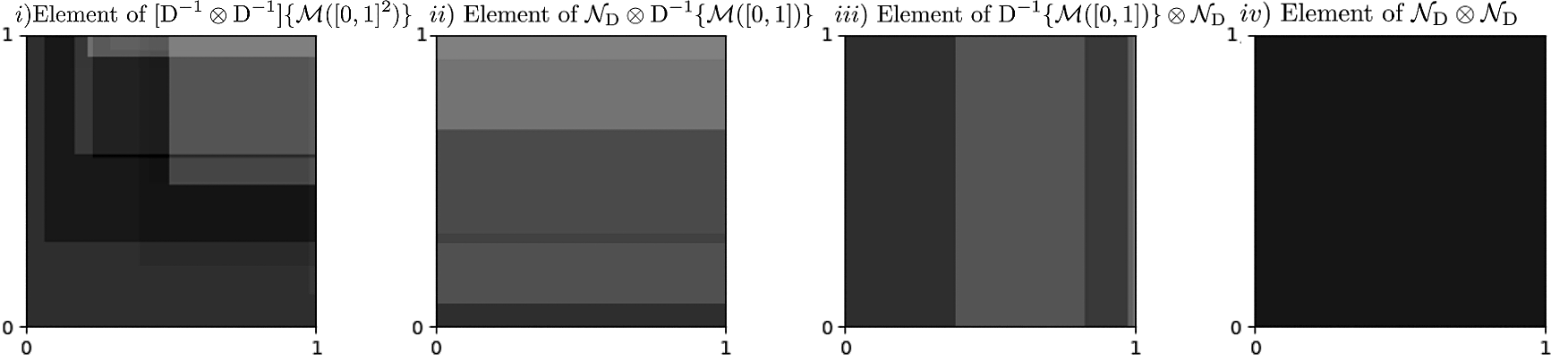}
    \caption{Direct-sum decomposition of the spline in Figure \ref{fig:1} into its four components, according to \eqref{eq:1.7}. The illustration is only made on the unit square $[0,1]^2$.}
    \label{fig:2}
\end{figure}

\subsection{Related Work}
\begin{itemize}

    \item \textbf{$\M$ Space Regularization.}
 Among the representer-theorem (RT) results \citep{boyer2019representer,bredies2020sparsity,flinth2019exact,unser2017splines} whose goals are to characterize the solution set $\mathcal{V}$, the results in \citep{bredies2020sparsity, flinth2019exact} are applicable to the regularization term $\norm{\Ls\{\cdot\}}_{\M(\R^2)}$, where $\Ls$ is an arbitrary partial-differential operator. There, the authors leave the null space $\mathcal{N}_{\Ls}$ regularization-free and conclude that there exists a solution that is an $\Ls$-spline with at most $(M-N)$ knots, where $N$ is the dimension of the vector space $\mathcal{A}=\{(p(\bb{t}_m))_{m=1}^M: p\in\mathcal{N}_{\Ls}\}$. This conceptually nice result has the benefit of having few constraints on the choice of $\Ls$. In turn, one limitation is that the absence of regularization on $\mathcal{N}_{\Ls}$ causes an ill-posedness of the IP, due to a potentially unbounded set of solutions and to the lack of coercivity of the loss functional.

    \item \textbf{Hilbert-Space Regularization.} A large part of the literature, originating from the work of Tikhonov\citep{tikhonov1963solution} and the RKHS theory \citep{scholkopf2001generalized,wahba1990spline}, formulates the regularization with an inner product. The author of \citep{unser2021unifying} showed that, in abstract Hilbert spaces and for a general class of admissible measurement operators, the solution is unique and can be described in terms of the Riesz map.
    In the particular case of pointwise evaluations, where the measurement functional is formed with Diracs, the solution becomes a linear combination of kernel elements and one falls back on the RKHS theory. Hence, when the kernel is defined in terms of tensor-product splines or B-splines \citep{scholkopf2002learning,vapnik1996support}, the RKHS theory advantageously promote such splines as solutions. Wendland in \citep{wendland2004scattered} develops this technique for the regularization with general partial-differential operators. Nonetheless, the knots there are fixed to the data localization, in contrast to the knots of our variational splines, which are free. In addition, the RKHS theory only permits, by nature, splines with even order $N_i$.
    
    \item \textbf{Vector-Valued Differential Operators.}
        In contrast to our use of scalar-valued PDO, one may choose a vector-valued differential operator for regularization. Total-variation (TV) \citep{chambolle2004algorithm,getreuer2012rudin,rudin1992nonlinear} and Hessian total-variation (HTV) \citep{aziznejad2023measuring,lefkimmiatis2013hessian} regularizations, associated to the gradient and Hessian operators, are the most salient examples. Because the vector-valued nature of these formulations automatically restrict the null space to be finite-dimensional, they are useful regularizations in dimensions higher than one. The solutions are typically non-separable, as studied in \citep{bredies2020sparsity} for TV, and in \citep{ambrosio2024linear} for HTV. This, in turn, leads to more difficult discretization and convergence analysis \citep{bartels2015total,pourya2024box}.
\end{itemize}

%%%%%%%%%%%%%%%%%%%%%%%%%%%%%%%%%%%%%%%%%%%%%%%%%%%%%
\section{Tensor-Product Background}
 For concision, we introduce tensor products in a sub-general fashion while, sometimes, slightly abusing the notation. A complete exposition is given in \citep{kubrusly2023bilinear,ryan2002introduction}.

\subsection{Notations and Vector Spaces}
\begin{itemize}
\item \emph{Indexing.} We shall use the variables $x_1$ and $x_2$ to index the domain of $f_1\otimes f_2$, the tensor product of two distributions. 
\item \emph{Natural Numbers.} The set of natural numbers, denoted by $\N$, is the set $\{1,2,\cdots\}$.
\item \emph{The Star.} The superscript $^{\star}$ is used to denote the adjoint of an operator, or the continuous dual space. 
\item \emph{Completion.} The symbol $\widehat{\otimes}$ is used to indicate the completion of a tensor product.
\item \emph{Distributions.} 
The space of compactly supported, infinitely smooth functions is denoted by  $\D(\R^D)$, while its continuous dual (the space of distribution) is $\D'(\R^D)$. The space of compactly supported distributions is denoted by $\mathcal{E}'(\R^D).$
\item \emph{Measures.} Let $\bb{K}\subset\R^D$ be a compact set. Then, the vector space $(\mathcal{M}(\bb{K}),\norm{\cdot}_{\M(\bb{K})})$ is the Banach space of bounded Radon measures supported on $\bb{K}$, equipped with the total-variation norm $\norm{\cdot}_{\M(\bb{K})}$ \citep{rudin}. In virtue of the Riesz-Markov-Kakutani theorem, it is the dual space of $\left(\C(\bb{K}),\norm{\cdot}_{\infty(\bb{K})}\right)$, the Banach space of functions that are continuous on $\bb{K}$, equipped with the sup norm. Likewise, $\M(\R^D)$ is the dual space of $\C_0(\R^D),$ the space of continuous functions that vanish at $\infty.$
\end{itemize}

\subsection{Tensor Product of Distributions}
\label{sec:2.2}
The tensor product $\otimes$ of two distributions $(f_1,f_2)\in\D'(\R)^2$ is the distribution $f_1\otimes f_2\in\D'(\R^2)$ such that 
\begin{equation}
\label{eq:2.1.1}
    \forall (\phi_1,\phi_2)\in\D(\R)^2:\quad\langle f_1\otimes f_2,\phi_1\otimes\phi_2\rangle=\langle f_1,\phi_1\rangle\langle f_2,\phi_2\rangle,
\end{equation}
where $\big(\phi_1\otimes\phi_2\big)(x_1,x_2)=\phi_1(x_1)\phi_2(x_2)$. The tensor-product space $\D'(\R)\otimes \D'(\R)\subset\D'(\R^2)$ is defined as 
\begin{equation}
    \D'(\R)\otimes \D'(\R):=\left\{\sum_{k=1}^Kf_{1,k}\otimes f_{2,k},\quad K\in\mathbb{N},(f_{1,k},f_{2,k})\in\D'(\R)^2\right\}.
\end{equation}
The tensor product of two infinite-dimensional vector spaces is never complete. The completion $\D'(\R)\widehat{\otimes}\D'(\R)$ (consistently, under the injective or the projective norm \citep{treves2016topological}) is identified to $\D'(\R^2).$ The tensor product $m_1\otimes m_2\in\M(\R^2)$ of two measures $(m_1,m_2)\in\M(\R)^2$, sometimes called the product measure \citep{halmos2013measure}, is defined likewise, with \eqref{eq:2.1.1} holding  $\forall(\phi_1,\phi_2)\in\C_0(\R)^2.$ The vector space $\D(\R)\otimes\D(\R)$ is defined by the multiplication of two functions, while its completion $\D(\R)\widehat{\otimes}\D(\R)$ under the injective or the projective norm is identified to $\D(\R^2)$ \citep{treves2016topological}.

\subsection{Tensor Product of Differential Operators}
\label{sec:2.3}
Let $\Ls_1:\D'(\R)\to\D'(\R)$ and $\Ls_2:\D'(\R)\to\D'(\R)$ be two ODO with null space $\mathcal{N}_{\Ls_1}$ and $\mathcal{N}_{\Ls_2}$. Then, $\Ls_1\otimes\Ls_2:\D'(\R^2)\to\D'(\R^2)$ is such that 
\begin{equation}
    \Ls_1\otimes\Ls_2:\begin{cases}
    \D'(\R)\otimes\D'(\R)\to\D'(\R^2)\\
    f=\sum_{k=1}^Kf_{1,k}\otimes f_{2,k}\mapsto[\Ls_1\otimes\Ls_2]\{f\}=\sum_{k=1}^K\Ls_1\{f_{1,k}\}\otimes\Ls_2\{f_{2,k}\},
    \end{cases}.
\end{equation}
Its action is extended to the whole vector space $\D'(\R^2)$ in reason of the density of $\D'(\R)\otimes\D'(\R)$. Further, $\Ls_1\otimes\Ls_2$ is naturally interpreted as a differential operator on $\R^2$, such that $\Ls_1$ differentiates in the $x_1$ dimension and $\Ls_2$ differentiates in the $x_2$ dimension. A distribution $g_{\Ls}\in\D'(\R)$ is a Green's function of the ODO $\Ls$ if $\Ls\{g_{\Ls}\}=\delta_0$. Likewise, a distribution $g_{\LL}\in\D'(\R^2)$ is a Green's function of the differential operator $\LL$ if 
\begin{equation}
    [\LL]\{g_{\LL}\}=\delta_0\otimes\delta_0=\delta_{\bb{0}},\quad\quad\bb{0}\in\R^2.
\end{equation}
Obverse that even though the causal Green's function is unique, a general Green's function is unique only up to an element of the null space $\mathcal{N}_{\LL}$ of $\LL$.
In Proposition \ref{prop:green} and Proposition \ref{prop:null space}, we characterize the causal Green's function and the null space of three differential operators of interest. As classic statements, they are provided without proof. 

\begin{Proposition} 
\label{prop:green}
Let $\Ls_1$ and $\Ls_2$ be two ODO with causal Green's function $g_{\Ls_1}$ and $g_{\Ls_2}$. Let $\I$ be the identity operator. Then,
     \begin{itemize}
        \item [1.]the function $g_{\Ls_1}\otimes g_{\Ls_2}$ is the unique causal Green's function of $\Ls_1\otimes\Ls_2;$
        \item [2.]the function $\delta_0\otimes g_{\Ls_2}$ is the unique causal Green's function of $\I\otimes\Ls_2;$
        \item [3.]the function $g_{\Ls_1}\otimes \delta_0$ is the unique causal Green's function of $\Ls_1\otimes\I.$
    \end{itemize}
\end{Proposition}

\begin{Proposition}
\label{prop:null space}
Let $\Ls_1$ and $\Ls_2$ be two ODO. Then,
     \begin{itemize}
        \item [1.]the null space $\mathcal{N}_{\Ls_1\otimes\Ls_2}$ of $\Ls_1\otimes\Ls_2$ is $\D'(\R)\otimes\mathcal{N}_{\Ls_2}+\mathcal{N}_{\Ls_1}\otimes\D'(\R)$;
        \item [2.]the null space $\mathcal{N}_{\I\otimes\Ls_2}$ of $\I\otimes\Ls_2$ is $\D'(\R)\otimes\mathcal{N}_{\Ls_2}$;
        \item [3.]the null space $\mathcal{N}_{\Ls_1\otimes\I}$ of $\Ls_1\otimes\I$ is $\mathcal{N}_{\Ls_1}\otimes\D'(\R)$.
    \end{itemize}
\end{Proposition}

\subsection{Tensor Product of Banach Spaces}
\label{sec:2.4}

Let $\mathcal{B}_1$ and $\mathcal{B}_2$ be two Banach subspaces of $\D'(\R)$, with norms $\norm{\cdot}_{\mathcal{B}_1}$ and $\norm{\cdot}_{\mathcal{B}_2}$. Define $\mathcal{B}_1^{\star},\mathcal{B}_2^{\star}$ as the dual spaces, equipped with the dual norm $\norm{\cdot}_{\mathcal{B}_1^{\star}},\norm{\cdot}_{\mathcal{B}_2^{\star}}$. Then, the tensor product $\mathcal{B}_1\otimes\mathcal{B}_2$ is defined as in \eqref{eq:2.1.1} and needs to be completed. Here, we make use of the injective norm \citep{ryan2002introduction}, noted $\norm{\cdot}_{\epsilon(\mathcal{B}_1,\mathcal{B}_2)}$ or, if context is clear, $\norm{\cdot}_{\epsilon}$, which is such that 
\begin{equation}
    \forall w\in \mathcal{B}_1\otimes\mathcal{B}_2:\quad\quad\norm{w}_{\epsilon}:=\underset{\underset{f_1\in\mathcal{B}_1^{\star}:\norm{f_1}_{\mathcal{B}_1^{\star}}=1}{f_2\in\mathcal{B}_2^{\star}:\norm{f_2}_{\mathcal{B}_2^{\star}}=1}}{\text{sup }}\left\vert\sum_{k=1}^K\langle f_1,v_{1,k}\rangle\langle f_2, v_{2,k}\rangle\right\vert,
\end{equation}
where $\sum_{k=1}^Kv_{1,p}\otimes v_{2,p}=w$ is a decomposition of $w$. It is such that, $\forall(v_1,v_2)\in\mathcal{B}_1\times\mathcal{B}_2:\quad$ $\norm{v_1\otimes v_2}_{\epsilon}=\norm{v_1}_{\mathcal{B}_1}\norm{v_2}_{\mathcal{B}_2}$ \citep[Proposition 3.1]{ryan2002introduction}.
The completion of $\mathcal{B}_1\otimes\mathcal{B}_2$ with respect to the injective norm is written $\mathcal{B}_1\widehat{\otimes}_{\epsilon}\mathcal{B}_2$ and is a Banach space, by definition. In Theorem \ref{th:densityC}, we recall two density results. The first one follows from \citep[page 49]{ryan2002introduction} and the second one is a simple extension, which is provided without proof.
\begin{Theorem}
    \label{th:densityC}
    The following completion results hold. 
\begin{itemize}
    \item [i.] For continuous functions with $\mathrm{T}=[-r,r]\subset\R,$ 
    \begin{equation}
        \left(\C(\mathrm{T})\widehat{\otimes}_\epsilon\C(\mathrm{T}),\norm{\cdot}_{\epsilon}\right)=\left(\C(\mathrm{T}^2),\norm{\cdot}_{\infty(\R^2)}\right).
    \end{equation}
    \item [ii.]For continuous functions that vanish at infinity,
    \begin{equation}
        \left(\C_0(\R)\widehat{\otimes}_\epsilon\C_0(\R),\norm{\cdot}_{\epsilon}\right)=\left(\C_0(\R^2),\norm{\cdot}_{\infty(\R^2})\right).
    \end{equation}
\end{itemize}
\end{Theorem}
In turn, the injective norm may not be adequate to complete the tensor product $\B_1^{\star}\otimes\B_2^{\star}$, in which case we make use of the $\w$ completion, noted $\B_1^{\star}\widehat{\otimes}_w\B_2^{\star}$. The latter is defined as the Banach space of all functionals $f$ for which there exists a sequence $(f_{\ell})_{\ell=1}^{\infty}$ with $f_{\ell}\in\B_1^{\star}\otimes\B_2^{\star}$ such that 
\begin{equation}
    \forall w\in\B_1\widehat{\otimes}_{\epsilon}\B_2,\quad\underset{\ell\to\infty}{\text{lim}}\langle f_{\ell},w\rangle=\langle f,w\rangle.
\end{equation}
Observe that the action $\langle f_{\ell},w\rangle$ of $f_{\ell}$ on $w$ is well-defined because
\begin{equation}
    \left(\B_1^{\star}\otimes\B_2^{\star},\norm{\cdot}_{\epsilon(\B_1,\B_2)^{\star}}\right)\subset\left(\B_1\widehat{\otimes}_{\epsilon}\B_2,\norm{\cdot}_{\epsilon(\B_1,\B_2)}\right)^{\star}.
\end{equation}
In Theorem \ref{th:3.1}, whose proof is given in \ref{app:2.3}, we apply this completion technique to the case of bounded measures.
\begin{Theorem}
\label{th:3.1}The tensor-product space of one-dimensional bounded measures is $\w$-dense in the space of two-dimensional bounded measures, with
\begin{equation}
    \left(\M(\R)\widehat{\otimes}_w\M(\R),\norm{\cdot}_{\epsilon^{\star}}\right)=\left(\M(\R^2),\norm{\cdot}_{\M(\R^2)}\right).
\end{equation}
\end{Theorem}

%%%%%%%%%%%%%%%%%%%%%%%%%%%%%%%%%%%%%%%%%%%%%%%%%%%%%
\section{Review of the One-Dimensional Theory}

\subsection{On the Real Line}
\label{sec:3}

As preparation for our two-dimensional tensor-product theory, we first review the one-dimensional theory \citep{unser2019native,unser2017splines}. The key results are synthesized, and some of them are slightly extended. Let $\Ls=(\Dd-\alpha\mathrm{I})^{N}$ be a differential operator with a single root $\alpha$ of multiplicity (or degree) $N$. The unique causal Green's function of $\Ls$ is $g_{\Ls}(\cdot)=\frac{(\cdot)_+^{N-1}e^{\alpha\cdot}}{(N-1)!}$ and its nullspace, of dimension $N$, is
\begin{equation}
    \mathcal{N}_{\Ls}=\{f\in\D'(\R):\Ls\{f\}=0\}=\mathrm{span}\left\{\frac{(\cdot)^{n-1}e^{\alpha \cdot}}{(n-1)!}\right\}_{n=1}^{N}.
\end{equation}
\begin{Definition}
\label{def:spline}
    A distribution $f\in\D'(\R)$ is an $\Ls$-spline if 
    \begin{equation}
        \Ls\{f\}=\sum_{k=1}^Ka_k\delta_{x_k},
    \end{equation}
    for some $K\in\mathbb{N}$, coefficients $(a_k)_{k=1}^K\in\R^K,$ and knots $(x_k)_{k=1}^K\in\R^K$. If we define the innovation $m=\sum_{k=1}^Ka_k\delta_{x_k}\in\D'(\R)$, we can generate such a spline as
    \begin{equation}
        f(\cdot)=(g_{\Ls}\ast m)(\cdot)+q(\cdot),\quad\text{with}\quad q\in\mathcal{N}_{\Ls}.
    \end{equation}
\end{Definition}
Splines live in the native space $\M_{\Ls}(\R)=\{f:\R\to\R,\quad\norm{\Ls\{f\}}_{\M}<\infty\}$, whose predual is denoted $\C_{\Ls}(\R)$. In order to review the properties of these spaces, it is necessary to introduce the notion of admissible systems $(\bb{p}, \bm{\phi})$, and their associated kernel $g_{\bm{\phi}}$. 
\begin{Definition}
\label{def:admissible}
    Let $\bb{p}=(p_n)_{n=1}^N$ be a basis of $\mathcal{N}_{\Ls}$ and $\bm{\phi}=(\phi_{n})_{n=1}^N$ be a complementary set of analysis functions. Then, $(\bb{p}, \bm{\phi})$ is a universal $\Ls$- admissible system if, $\forall n\in[1\ldots N]$, the function $\phi_n\in\mathcal{D}(\R)$ is compactly supported with $\mathrm{supp}(\phi_n)\subset[\phi^-,\phi^+]$ and is such that 
    \begin{equation}
    \label{eq:3.13}
            \forall (n,m)\in[1,N]^2:\quad\langle 
            p_n, \phi_m\rangle=\delta_{n-m}.
    \end{equation}
\end{Definition}
Observe that $\bm{\phi}$ is a basis for the dual space $\mathcal{N}_{\Ls}^{\star}$ of $\mathcal{N}_{\Ls}$. The projection operators $\text{proj}_{\mathcal{N}_{\Ls}^{\star}}$ and $\text{proj}_{\mathcal{N}_{\Ls}}$ are defined as 
\begin{equation}
    \text{proj}_{\mathcal{N}_{\Ls}}:
    \begin{cases}
        \D'(\R)\to\mathcal{N}_{\Ls}\\
        f\mapsto\sum_{n=1}^N\langle f,\phi_n\rangle p_n,
    \end{cases}\quad\quad    \text{proj}_{\mathcal{N}_{\Ls}^{\star}}:
    \begin{cases}
        \D(\R)\to\mathcal{N}_{\Ls}^{\star}\\
        f\mapsto\sum_{n=1}^N\langle p_n,f\rangle\phi_n.
    \end{cases}
\end{equation}
 The kernel $g_{\bm{\phi}}$ is defined as (see \citep{unser2017splines}) 
\begin{equation}
\label{eq:2.8}
  g_{\bm{\phi}}(t,x)=g_{\Ls}(t-x)-\sum_{n=1}^N(g_{\Ls}^{\vee}\ast\phi_n)(x)p_n(t)=g_{\Ls}(t-x)-\sum_{n=1}^N\langle g_{\Ls}(\cdot),\phi_{n}(\cdot+x)\rangle p_n(t), 
\end{equation}
and the associated integral transforms $\Ls_{\bm{\phi}}^{-1}$ and $\Ls_{\bm{\phi}}^{-\star}$ are 
\begin{equation}
\label{eq:2.9}
     \Ls_{\bm{\phi}}^{-1}\{f\}=\int_{\R}g_{\bm{\phi}}(t,x)f(x)\mathrm{d}x,\quad\quad \Ls_{\bm{\phi}}^{-\star}\{g\}=\int_{\R}g_{\bm{\phi}}(t,x)g(t)\mathrm{d}t,
\end{equation}
under the implicit assumption of well-posedness. We recall Theorem \ref{th:1} and \ref{Th:predual1D} that establish several key properties of $\M_\Ls(\R)$ and $\C_{\Ls}(\R)$. The vector spaces $\mathcal{N}_{\Ls},\mathcal{N}_{\Ls}^{\star}$ are equipped with the norms $\norm{\cdot}_{\bb{p}}$,$\norm{\cdot}_{\bm{\phi}}$ defined as
\begin{equation}
    \norm{\sum_{n=1}^Na_np_n}_{\bb{p}}:=\sum_{n=1}^{N}\vert a_n\vert,\quad\quad\norm{\sum_{n=1}^Na_n\phi_n}_{\bm{\phi}}:=\underset{1\leq n\leq N}{\text{max}}\vert a_n\vert.
\end{equation}
We recall that two vector spaces $\mathcal{B}_1$ and $\mathcal{B}_1$ are in direct sum if the map 
\begin{equation}
    \mathrm{S}:\mathcal{B}_1\times\mathcal{B}_2\to\mathcal{B}_1+\mathcal{B}_2,\quad(v_1,v_2)\mapsto v_1+v_2
\end{equation}
is an isomorphism.
\begin{Theorem}[\citep{unser2017splines}]
\label{th:1}
If $(\bb{p}, \bm{\phi})$ is a universal $\Ls$-admissible system, then the space $\M_{\Ls}(\R)$ admits the direct-sum decomposition
\begin{equation}
\label{eq:3.23}
    \M_{\Ls}(\R)=\Ls_{\bm{\phi}}^{-1}\{\M(\R)\}\oplus\mathcal{N}_{\Ls},\quad\text{with the norm}\quad \norm{\cdot}_{\M_{\Ls}}:=\norm{\Ls\{\cdot\}}_{\M}+\norm{\mathrm{proj}_{\mathcal{N}_{\Ls}}\{\cdot\}}_{\bb{p}},
\end{equation}
for which it is a Banach space. In particular, it holds that
\begin{equation}
     \label{eq:3.38}
        \forall m\in\M(\R):\quad\mathrm{proj}_{\mathcal{N}_{\Ls}}\{\Ls_{\bm{\phi}}^{-1}\{m\}\}=0.
\end{equation}
Furthermore, the extreme points of the centred unit ball in $(\M_{\Ls}(\R), \norm{\Ls\{\cdot\}}_{\M})$ are exactly $\Ls$-splines with one knot $x_1$ and amplitude $a_1\in\{-1,1\}$.
\end{Theorem}

\begin{Theorem}[\citep{unser2017splines}]
\label{Th:predual1D}
If $(\bb{p}, \bm{\phi})$ is a universal admissible system, then the function space
\begin{equation}
\label{eq:3.2.40}
    \C_{\Ls}(\R)=\left\{g=\Ls^{\star}\{v\}+\sum_{n=1}^Na_n\phi_n:v\in\C_{0}(\R),\bb{a}=(a_n)_{n=1}^N\in\R^N\right\}
\end{equation}
is such that every $g\in\C_{\Ls}(\R)$ has a unique direct-sum representation as in \eqref{eq:3.2.40}; in particular,    \begin{center}
    \vspace*{-\baselineskip}
    \begin{minipage}{0.5\linewidth}
    \begin{align}
            \Ls_{\bm{\phi}}^{-\star}\left\{\sum_{n=1}^Na_n\phi_n\right\}=0,\label{eq:3.43}\\
            \Ls_{\bm{\phi}}^{-\star}\{\Ls^{\star}\{v\}\}=v\label{eq:3.43bis},
    \end{align}
    \end{minipage}\begin{minipage}{0.5\linewidth}
        \begin{align}
            \mathrm{proj}_{\mathcal{N}_{\Ls}^{\star}}\{\Ls_{\bm{\phi}}^{-\star}\{v\}\}=0\label{eq:3.44}.
        \end{align}
        \end{minipage}
    \end{center}
Furthermore, the space $\C_{\Ls}(\R)$ is a Banach space. It is also the predual of $\M_{\Ls}(\R)$ for the norm
    \begin{equation}
    \norm{g}_{\C_{\Ls}}:=\mathrm{max}\left(\norm{\Ls_{\bm{\phi}}^{-\star}\{g\}}_{\infty(\R)},\norm{\mathrm{proj}_{\mathcal{N}_{\Ls}^{\star}}\{g\}}_{\bm{\phi}}\right)=\mathrm{max}(\norm{v}_{\infty(\R)},\norm{\bb{a}}_{\infty}).
    \end{equation}
\end{Theorem}

This characterization of the predual $\C_{\Ls}(\R)$ enables us to extend the notion of admissible system to a collection $(\phi_n)_{n=1}^N$ with $\phi_{n}\in\C_{\Ls}(\R)$ \citep{unser2019native}. Indeed, we observe \emph{a posteriori} that the construction of the spaces $\M_{\Ls}(\R)$ and $\C_{\Ls}(\R),$ as well as Theorem \ref{th:1} and Theorem \ref{Th:predual1D}, hold for any $\Ls$-admissible system $(\bb{p}, \bm{\phi})$. We conclude this section with Proposition \ref{prop:compactsupport1D} that describes the regularity and the support of $g_{\bm{\phi}}$. The support part is proved in \citep{guillemet2025adaptive}, and the regularity part can directly be inferred by looking at \eqref{eq:2.8}.
\begin{Proposition}
\label{prop:compactsupport1D}
If $(\bb{p}, \bm{\phi})$ is an $\Ls$-admissible system, then the regularity of the kernel $g_{\bm{\phi}}$ is such that 
\begin{itemize}
    \item for all $t\in\R$, $g_{\bm{\phi}}(t,\cdot)\in\mathcal{L}_{\infty,c}(\R)$ and $\langle g_{\Ls}(\cdot),\phi_n(\cdot+x)\rangle$ is continuous in $x$;
    \item for all $x\in\R$, $\frac{\partial^{N-1}g_{\bm{\phi}}(\cdot,x)}{\partial t^{N-1}}\in\mathcal{L}_{\infty,\mathrm{loc}}(\R)$ and, if $N>1$, $g_{\bm{\phi}}(\cdot,x)\in\C^{N-2}(\R)$.
\end{itemize}
    In addition, the support of the kernel $g_{\bm{\phi}}$ is such that 
   \begin{equation}
       \forall t\in\R,\forall x\notin[\mathrm{min}(t,\phi^-),\mathrm{max}(t,\phi^+)]:\quad g_{\bm{\phi}}(t,x)=0
   \end{equation}
   or, equivalently,
   \begin{equation}
       \begin{cases}
           \forall x<\phi^-,\forall t>x\\
           \forall x>\phi^+,\forall t<x
        \end{cases}:\quad g_{\bm{\phi}}(t,x)=0.
   \end{equation}
\end{Proposition}

\subsection{On an Interval}
\label{sec:fund}

Formulating the theory directly on the space $\M_{\Ls}(\R)$ is not always convenient because of the unbounded domain and underlying non-shift-invariant kernel $g_{\bm{\phi}}$. Instead, one can optimise over the localised search space 
\citep{guillemet2025convergence}: 
\begin{equation}
    \M_{\Ls}(\K)=\{f\in\mathcal{D}'(\R):\quad\Ls\{f\}\in\M(\K)\},
\end{equation}
where $\K=[K^-,K^+]$. In this case, one can specify, under meaningful constraints, the notion of an admissible system. Localised and fundamental systems were introduced and studied in \citep{guillemet2025sampling} for $\Ls=\Dd^N$. We extend these to the case $\Ls=(\Dd-\alpha\I)^N$.

First, localised systems $(\bb{p},\bm{\phi})$ yield an optimal weight distribution of the kernel, in the sense of Proposition \ref{prop:localizedjustification}.

\begin{Definition}
\label{def:localizedsystem}
    An $\Ls$-admissible system $(\bb{p},\bm{\phi})$ is called $\K$-localised if $[\phi^-,\phi^+]\subset[K^-,K^+]$.
\end{Definition}

\begin{Proposition}[\citep{guillemet2025sampling}]
\label{prop:localizedjustification}
    If $(\bb{p},\bm{\phi})$ is an $\Ls$-admissible system, then 
\begin{equation}
\label{eq:loc.1}
    \forall x\notin[\mathrm{min}(\phi^-,K^-),\mathrm{max}(\phi^+,K^+)],\forall t\in[K^-,K^+]:\quad g_{\bm{\phi}}(t,x)=0.
\end{equation}
    If $(\bb{p},\bm{\phi})$ is a $\K$-localised system, then 
\begin{equation}
\label{eq:loc.2}
    \forall x\notin[K^-,K^+],\forall t\in[K^-,K^+]:\quad g_{\bm{\phi}}(t,x)=0.
\end{equation}
\end{Proposition}

Second, $\K$-fundamental systems $(\bb{p},\bm{\iota})$ yield a kernel $g_{\bm{\iota}}$ that is locally shift-invariant over $\K$.

\begin{Definition}
\label{def:idealsystem}
    A $\K$-localised system $(\bb{p},\bm{\phi})$ is called $\K$-fundamental if
    \begin{itemize}
        \item [1.] for all $n\in[1\ldots N]$, the functional $\phi_n:\M_{\Ls}(\K)\to\R$ is $\w$-continuous;
        \item [2.] the kernel $g_{\bm{\phi}}$ reduces to the Green's function, in the sense that
    \begin{equation}
    \forall x\in[K^-,\infty[:\quad g_{\bm{\phi}}(\cdot,x)=g_{\Ls}(\cdot-x);
    \end{equation}
        \item [3.] the null-space element $p_{n}$ is the shifted canonical exponential polynomial, namely
    \begin{equation}
    \label{eq:4.4.71}
        p_{n}(\cdot)=\frac{(\cdot-K^-)^{n-1}\mathrm{e}^{\alpha(\cdot-K^-)}}{(n-1)!}.
    \end{equation}
    \end{itemize}
\end{Definition}

To describe the unique $\K$-fundamental system, we recall the definition of the left-limit derivative $\Dd^{N-1,\star}\{\delta_t^-\}$, 
\begin{equation}
\label{eq:2.6}
\langle\Dd^{N-1,\star}\{\delta_{t}^-\},\cdot\rangle:\begin{cases}
\M_{\Ls}(\K)\mapsto\R,\\
f\mapsto\underset{\epsilon\downarrow0}{\mathrm{lim }}\,f^{(N-1)}(t-\epsilon).
\end{cases}
\end{equation}

\begin{Theorem}[\citep{guillemet2025sampling}]
 \label{prop:ideal2}
The system $(\bb{p},\bm{\iota})$ defined as follows is the unique $\K$-fundamental system: the null-space element $p_{i,n}$ is given by \eqref{eq:4.4.71}, and $\iota_{i,n}$ is such that, $\forall f\in\M_{\Ls}(\K)$,
\begin{align}
   \forall n\in[1\ldots N-1],\quad &\langle f,\iota_{n}\rangle=\left\langle\delta_{K^-},\left(\Dd-\alpha\mathrm{I}\right)^{n-1}\{f\}
    \right\rangle,\nonumber\\
    &\langle f,\iota_{N}\rangle=\left\langle\delta_{K^-}^-,\left(\Dd-\alpha\mathrm{I}\right)^{N-1}\{f\}
    \right\rangle,
\end{align}
where $\delta_{K^-}^-$ denotes the left limit at $K^-$.
\end{Theorem}

A consequence of the local shift-invariance is that, in contrast with functions in $\M_{\Ls}(\R)$, functions in $\M_{\Ls}(\K)$ admit a direct-sum, integral representation with a shift-invariant kernel. 

\begin{Proposition}[\citep{guillemet2025sampling}]
If $(\bb{p},\bm{\iota})$ is the $\K$-fundamental system, then the Banach space $\M_{\Ls}(\K)$ satisfies 
\begin{equation}
\label{eq:fund.36}
    \M_{\Ls}(\K)=\Ls_{\bm{\iota}}^{-1}\{\M(\K)\}\oplus\mathcal{N}_{\Ls}=\Ls^{-1}\{\M(\K)\}\oplus\mathcal{N}_{\Ls},
\end{equation}
where $\Ls^{-1}$ is the integral operator whose kernel is $g_{\Ls}$ and 
\begin{equation}
    \Ls^{-1}\{\M(\K)\}=\left\{\int_{\K}g_{\Ls}(\cdot-x)\mathrm{d}m(x): m\in\M(\K)\right\}\label{eq:4.4.76}.
\end{equation}
\end{Proposition}

%%%%%%%%%%%%%%%%%%%%%%%%%%%%%%%%%%%%%%%%%%%%%%%%%%%%%
\section{Two-Dimensional Theory}
\label{sec:4}
We now move on to the tensor product of two differential operators $(\Dd-\alpha_i\mathrm{I})^{N_i}$, for $i\in\{1,2\}$, with a nullspace $\mathcal{N}_{\Ls_i}$ of dimension $N_i$, spanned by a basis $\bb{p}_i=(p_{i,n})_{n=1}^{N_i}$. The construction of $\M_{\LL}(\R^2)$ rests on our ability to ``smoothen", or integrate two-dimensional measures with respect to $[\LL]$. Therefore, in order to synthesize a function from a Dirac mass, we introduce Definition \ref{def:universalsystem} and extend the technique used in \citep{unser2017splines}.

\begin{Definition}
\label{def:universalsystem}
    A system  $((\bb{p}_1,\bb{p}_2),(\bm{\phi}_1,\bm{\phi}_2))$ is a universal $[\LL]$-admissible system if, for $i\in\{1,2\}$, $(\bb{p}_i,\bm{\phi}_i)$  is a universal $\Ls_i$-admissible system.
\end{Definition}

Each universal $[\LL]$-admissible system has three associated kernels $g_{\bm{\phi}_1},g_{\bm{\phi}_2},$ and $g_{\bm{\phi}_1}\otimes g_{\bm{\phi}_2}$ as in \eqref{eq:2.8}, whose associated integral transforms are noted $\Ls_{\bm{\phi}_1}^{-1},\Ls_{\bm{\phi}_2}^{-1},$ and $\Ls_{\bm{\phi_1}}^{-1}\otimes\Ls_{\bm{\phi}_2}^{-1}$. These transforms are right-inverses of their associated differential operators, in the sense that
\begin{align}
&\Ls_1\circ \Ls^{-1}_{\bm{\phi_1}}=\mathrm{Id},\quad\text{and}\quad\Ls_2\circ\Ls^{-1}_{\bm{\phi_2}}=\mathrm{Id}\quad\text{over}\quad\M(\R),\label{eq:4.0.47}\\
&[\LL]\circ\left[\Ls^{-1}_{\bm{\phi_1}}\otimes\Ls^{-1}_{\bm{\phi_2}}\right]=\mathrm{Id}\quad\text{over}\quad\M(\R^2)\label{eq:4.0.48}.
\end{align}
Equation \eqref{eq:4.0.47} follows from Theorem \ref{th:1}, and \eqref{eq:4.0.48} follows from Theorem \ref{th:native} below.

\subsection{Predual}
\label{sec:4.1}
We now define the predual
\begin{align}
\label{eq:4.1.35}
    \C_{\Ls_1\otimes\Ls_2}(\R^2) = &\underbrace{\left(\left[\Ls_1^{\star}\otimes\Ls_2^{\star}\right]\left\{\C_0(\R^2)\right\}\right)}_{(i)}+\underbrace{\left(\mathcal{N}_{\Ls_1}^{\star}\otimes\Ls_2^{\star}\{\C_0(\R)\}\right)}_{(ii)}\nonumber\\
    &+\underbrace{\left(\Ls_1^{\star}\left\{\C_0(\R)\right\}\otimes\mathcal{N}_{\Ls_2}^{\star}\right)}_{(iii)}+\underbrace{\left(\mathcal{N}_{\Ls_1}^{\star}\otimes\mathcal{N}_{\Ls_2}^{\star}\right)}_{(iv)}.
\end{align}
In Proposition \ref{prop:isom}, whose proof is in \ref{app:4.1}, we study and norm the sub-component spaces of \eqref{eq:4.1.35}.
\begin{Proposition}
\label{prop:isom}\hspace{+1cm}
\begin{itemize}
    \item [1.] The vector space $(i)$  expressed as $\big\{\left[\Ls_1^{\star}\otimes\Ls_2^{\star}\right]\{v\}:\quad v\in\C_0(\R^2)\big\}$ and equipped with the norm $\norm{\left[\Ls_{\bm{\phi}_1}^{-\star}\otimes \Ls_{\bm{\phi}_2}^{-\star}\right]\{\cdot\}}_{\infty(\R^2)}$  is a Banach space. It is isometric to $\left(\C_0(\R^2),\norm{\cdot}_{\infty(\R^2)}\right)$ under the mapping $\Ls_{\bm{\phi}_1}^{-\star}\otimes \Ls_{\bm{\phi}_2}^{-\star}$ whose inverse is $\Ls_1^{\star}\otimes\Ls_2^{\star}.$ Furthermore, $(i)$ is the injective norm completion of $\Ls_1^{\star}\{\C_0(\R)\}\otimes\Ls_2^{\star}\{\C_0(\R)\}$.

    \item [2.]
The vector space  $(ii)$ expressed as $\left\{\sum_{n=1}^{N_1}\phi_{1,n}\otimes\Ls_{2}^{\star}\{v_{1,n}\}:\quad(v_{1,n})_{n=1}^{N_1}\in\C_0(\R)^{N_1}\right\}$ and equipped with the norm 
    \begin{align}
        &\norm{\left[\mathrm{proj}_{\mathcal{N}_{\Ls_1}^{\star}}\otimes \Ls_{\bm{\phi}_2}^{-\star}\right]\left\{\sum_{n=1}^{N_1}\phi_{1,n}\otimes\Ls_{2}^{\star}\{v_{1,n}\}\right\}}_{\bm{\phi}_1,\infty(\R)}\nonumber\\
        =&\norm{\sum_{n=1}^{N_1}\phi_{1,n}\otimes v_{1,n}}_{\bm{\phi}_1,\infty(\R)}:=\underset{x\in\R}{\mathrm{sup }}\underset{1\leq n\leq N_1}{\mathrm{max }}\vert v_{1,n}(x)\vert
    \end{align}
    is a Banach space. It is isometric to $\left(\mathcal{N}_{\Ls_1}^{\star}\otimes \C_0(\R),\norm{\cdot}_{\bm{\phi}_1,\infty(\R)}\right)$ under the mapping $\mathrm{proj}_{\mathcal{N}_{\Ls_1}^{\star}}\otimes \Ls_{\bm{\phi}_2}^{-\star}$ whose inverse is $\mathrm{proj}_{\mathcal{N}_{\Ls_1}^{\star}}\otimes\Ls_2^{\star}$.
        
    \item [3.] The vector space $(iii)$  expressed as $\left\{\sum_{n=1}^{N_2}\Ls_1^{\star}\{v_{2,n}\}\otimes\phi_{2,n}:\quad(v_{2,n})_{n=1}^{N_2}\in\C_0(\R)^{N_2}\right\}$ and equipped with the norm 
    \begin{align}
        &\norm{\left[\Ls_{\bm{\phi}_1}^{-\star}\otimes\mathrm{proj}_{\mathcal{N}_{\Ls_2}^{\star}}\right]\left\{\sum_{n=1}^{N_2}\Ls_{1}^{\star}\{v_{2,n}\}\otimes\phi_{2,n}\right\}}_{\infty(\R),\bm{\phi}_2}\nonumber\\
        =&\norm{\sum_{n=1}^{N_2} v_{2,n}\otimes\phi_{2,n}}_{\infty(\R),\bm{\phi_2}}:=\underset{1\leq n\leq N_2}{\mathrm{max }}\underset{x\in\R}{\mathrm{sup }}\vert v_{2,n}(x)\vert
    \end{align}
    is a Banach space. It is isometric to $\left(\C_0(\R)\otimes\mathcal{N}_{\Ls_2}^{\star},\norm{\cdot}_{\infty(\R),\bm{\phi}_2}\right)$ under the mapping $\Ls_{\bm{\phi}_1}^{-\star}\otimes\mathrm{proj}_{\mathcal{N}_{\Ls_2}^{\star}}$ whose inverse is $\Ls_1^{\star}\otimes\mathrm{proj}_{\mathcal{N}_{\Ls_2}^{\star}}$. 
    
    \item [4.] The vector space $(iv)$ expressed as
    \begin{equation}
        \left\{\sum_{n,n'=1}^{N_1,N_2}d_{n,n'}\phi_{1,n}\otimes\phi_{2,n'}:\quad(d_{n,n'})_{n,n'=1,1}^{N_1,N_2}\in\R^{N_1\times N_2}\right\}
    \end{equation}
     and equipped with the norm
    \begin{equation}
        \norm{\sum_{n,n'=1}^{N_1,N_2}d_{n,n'}\phi_{1,n}\otimes\phi_{2,n'}}_{\bm{\phi}_1,\bm{\phi}_2}:=\underset{1\leq n\leq N_1;1\leq n'\leq N_2}{\mathrm{max}}\vert d_{n,n'}\vert
    \end{equation}
    is a Banach space.
\end{itemize}
\end{Proposition}

In Theorem \ref{th:predual}, whose proof is in \ref{app:4.1}, we concatenate the results of Proposition \ref{prop:isom} to provide $\C_{\LL}(\R^2)$ with a direct-sum structure, which will prove helpful to characterize admissible measurement functionals. 
\begin{Theorem}
    \label{th:predual}
     The spaces $(i),(ii),(iii),$ and $(iv)$ are in direct sum and every $g\in \C_{\LL}(\R^2)$ has a unique direct-sum representation as \begin{equation}
        g=[\LL]\{v\}+\sum_{n=1}^{N_1}\phi_{1,n}\otimes \Ls_2^{\star}\{v_{1,n}\}+\sum_{n=1}^{N_2}\Ls_1^{\star}\{v_{2,n}\}\otimes\phi_{2,n}+\sum_{n,n'=1}^{N_1,N_2}d_{n,n'}\phi_{1,n}\otimes\phi_{2,n},\label{eq:4.1.43}
    \end{equation}
    for some $v\in \C_0(\R^2),$  $(v_{1,n})_{n=1}^{N_1}\in\C_0(\R)^{N_1},$ $(v_{2,n})_{n=1}^{N_2}\in\C_0(\R)^{N_2},$ and $(d_{n,n'})_{n,n'=1,1}^{N_1,N_2}\in\R^{N_1\times N_2}$. In addition, the mappings
    \begin{equation}
        \Ls_{\bm{\phi}_1}^{-\star}\otimes \Ls_{\bm{\phi}_2}^{-\star},\quad \mathrm{proj}_{\mathcal{N}_{\Ls_1}^{\star}}\otimes \Ls_{\bm{\phi}_2}^{-\star},\quad\Ls_{\bm{\phi}_1}^{-\star}\otimes\mathrm{proj}_{\mathcal{N}_{\Ls_2}^{\star}},\quad\mathrm{proj}_{\mathcal{N}_{\Ls_1}^{\star}}\otimes\mathrm{proj}_{\mathcal{N}_{\Ls_2}^{\star}}
    \end{equation}
    extend continuously from their domain of definition to the whole predual $\C_{\LL}(\R^2).$ Consequently, $\C_{\Ls_1\otimes\Ls_2}(\R^2)$ is a Banach space for the norm
    \begin{align}
        \norm{g}_{\C_{\Ls_1\otimes\Ls_2}}&=\mathrm{max}\bigg(\norm{\left[\Ls_{\bm{\phi}_1}^{-\star}\otimes \Ls_{\bm{\phi}_2}^{-\star}\right]\{g\}}_{\infty(\R^2)},\norm{\left[\mathrm{proj}_{\mathcal{N}_{\Ls_1}^{\star}}\otimes \Ls_{\bm{\phi}_2}^{-\star}\right]\{g\}}_{\bm{\phi}_1,\infty(\R)},\nonumber\\
        &\hspace{+1.4cm}\norm{\left[\Ls_{\bm{\phi}_1}^{-\star}\otimes\mathrm{proj}_{\mathcal{N}_{\Ls_2}^{\star}}\right]\{g\}}_{\infty(\R),\bm{\phi}_2},\norm{ \left[\mathrm{proj}_{\mathcal{N}_{\Ls_1}^{\star}}\otimes\mathrm{proj}_{\mathcal{N}_{\Ls_2}^{\star}}\right]\{g\}}_{\bm{\phi}_1,\bm{\phi}_2}\bigg)\nonumber\\
        &=\mathrm{max}\left(\norm{v}_{\infty(\R^2)},\underset{x\in\R,1\leq n\leq N_1}{\mathrm{sup }}\vert v_{1,n}(x)\vert,\underset{x\in\R,1\leq n\leq N_2}{\mathrm{sup }}\vert v_{2,n}(x)\vert,\norm{\bb{d}}_{\infty}\right)\label{eq:4.1.44}.
    \end{align}
\end{Theorem} 

From the proof of Theorem \ref{th:predual}, we also know that $\C_{\LL}(\R^2)$ is the $\norm{\cdot}_{\C_{\LL}}$-completion of $\C_{\Ls_1}(\R)\otimes\C_{\Ls_2}(\R)$, which emphasizes its tensor-product structure. In addition, $\C_{\LL}(\R^2)$ is, as a set, independant of the choice of the universal $[\LL]$-admissible system. In Corollary \ref{coro:admissible}, we provide a simple-to-verify condition to assess whether a compactly supported distribution is in $\C_{\LL}(\R)$, or not.
\begin{Corollary}
\label{coro:admissible}
    If $g\in\mathcal
    E'(\R^2)$, then $g\in\C_{\Ls_1\otimes\Ls_2}(\R^2)$ if and only if the following three conditions are verified
    \begin{itemize}
        \item [1.] the distribution  $\left[\Ls_{\bm{\phi}_1}^{-\star}\otimes \Ls_{\bm{\phi}_2}^{-\star}\right]\{g\}\in\mathcal{E}'(\R^2)$ is in $\C(\R^2);$
        \item [2.] the distribution $\left[\mathrm{proj}_{\mathcal{N}_{\Ls_1}^{\star}}\otimes \Ls_{\bm{\phi}_2}^{-\star}\right]\{g\}\in\mathcal{E}'(\R^2)$ is in $\C(\R^2);$
        \item [3.] the distribution $\left[\Ls_{\bm{\phi}_1}^{-\star}\otimes\mathrm{proj}_{\mathcal{N}_{\Ls_2}^{\star}}\right]\{g\}\in\mathcal{E}'(\R^2)$ is in $\C(\R^2).$
    \end{itemize}
Moreover, if $g=g_1\otimes g_2$ with $(g_1,g_2)\in\mathcal{E}'(\R)^2$, then $g\in\C_{\Ls_1\otimes\Ls_2}(\R^2)$ if and only if the following two conditions are verified
\begin{itemize}
    \item [4.] the distribution $g_{\Ls_1}^{\vee}\ast g_1\in\D'(\R)$ is in $\C(\R)$;
    \item [5.] the distribution $g_{\Ls_2}^{\vee}\ast g_2\in\D'(\R)$ is in $\C(\R)$.
\end{itemize}
\end{Corollary}
This shows that it is sufficient to verify the inclusion in $\C(\R^2)$, instead of $\C_0(\R^2).$
In Corollary \ref{coro:testfunctions}, we build on Corollary \ref{coro:admissible} to show that $\C_{\LL}(\R^2)$ is sufficiently large to include all Lebesgue functions, as well as evaluation functionals $\delta_{\bb{x}}$, provided that $\Ls_1$ and $\Ls_2$ are of degree greater than 1. Its proof is given in \ref{app:4.1}.

\begin{Corollary}
\label{coro:testfunctions}
The following inclusions hold
\begin{itemize}
    \item [1.] for infinitely smooth functions, $\D(\R^2)\subset\C_{\LL}(\R^2);$
    \item [2.] for integrable functions, $\forall p\in[1,\infty],\quad\mathcal{L}_{p,c}(\R^2)\subset\C_{\LL}(\R^2).$  
\end{itemize}       
If, in addition, the orders $N_1,N_2$ of $\Ls_1,\Ls_2$ are both greater than $1$, then 
\begin{itemize}
    \item [3.] for any $\bb{x}\in\R^2$, $\delta_{\bb{x}}\in\C_{\LL}(\R^2).$ 
\end{itemize}
\end{Corollary}

\subsection{Native Space}
\label{sec:4.2}
We now define the native space 
\begin{align}
    \M_{\Ls_1\otimes\Ls_2}(\R^2)=&\underbrace{\left(\left[\Ls_{\bm{\phi}_1}^{-1}\otimes \Ls_{\bm{\phi}_2}^{-1}\right]\{\M(\R^2)\}\right)}_{(i)^{\star}}+\underbrace{\left(\mathcal{N}_{\Ls_1}\otimes \Ls_{\bm{\phi}_2}^{-1}\{\M(\R)\}\right)}_{(ii)^{\star}}\nonumber\\
    &+\underbrace{\left(\Ls_{\bm{\phi}_1}^{-1}\{\M(\R)\}\otimes\mathcal{N}_{\Ls_2}\right)}_{(iii)^{\star}}+\underbrace{\left(\mathcal{N}_{\Ls_1}\otimes\mathcal{N}_{\Ls_2}\right)}_{(iv)^{\star}}\label{eq:4.2.69}.
\end{align}
Proposition \ref{prop:isom2}, whose proof is in \ref{app:4.2}, gives the structure of the sub-component spaces of \eqref{eq:4.2.69}.
\begin{Proposition}\hspace{+1cm}
\label{prop:isom2}
    \begin{itemize}
        \item [1.] The vector space $        (i)^{\star}$ expressed as $\left\{\left[\Ls_{\bm{\phi}_1}^{-1}\otimes \Ls_{\bm{\phi}_2}^{-1}\right]\{m\},\quad m\in\M(\R)^2\right\},$ where
        \begin{equation}
        \label{eq:4.2.47}
            \left[\Ls_{\bm{\phi}_1}^{-1}\otimes \Ls_{\bm{\phi}_2}^{-1}\right]\{m\}(\cdot)=\int_{\R^2}g_{\bm{\phi}_1}(\cdot,x_1)g_{\bm{\phi}_2}(\cdot,x_2)\mathrm{d}m(x_1,x_2),
        \end{equation}
        and equipped with the norm
    \begin{equation}
        \norm{[\Ls_1\otimes\Ls_2]\left\{\left[\Ls_{\bm{\phi}_1}^{-1}\otimes \Ls_{\bm{\phi}_2}^{-1}\right]\{m\}\right\}}_{\M(\R^2)}=\norm{m}_{\M(\R^2)},
    \end{equation}
    is a Banach space. It is isometric to $\left(\M(\R^2),\norm{\cdot}_{\M(\R^2)}\right)$, under the mapping $\LL$ whose inverse is  $\Ls_{\bm{\phi}_1}^{-1}\otimes \Ls_{\bm{\phi}_2}^{-1}$. 
    
    \item [2.] The vector space $(ii)^{\star}$ expressed as $\left\{\sum_{n=1}^{N_1}p_{1,n}\otimes \Ls_{\bm{\phi}_2}^{-1}\{m_{1,n}\},\quad(m_{1,n})_{n=1}^{N_1}\in\M(\R)^{N_1}\right\}$ and equipped with the norm
    \begin{align}
        &\norm{\left[\mathrm{proj}_{\mathcal{N}_{\Ls_1}}\otimes\Ls_2\right]\left\{\sum_{n=1}^{N_1}p_{1,n}\otimes \Ls_{\bm{\phi}_2}^{-1}\{m_{1,n}\}\right\}}_{\bb{p}_1,\M(\R)}\nonumber\\
        =&\norm{\sum_{n=1}^{N_1}p_{1,n} \otimes m_{1,n}}_{\bb{p}_1,\M(\R)}:=\sum_{n=1}^{N_1}\norm{m_{1,n}}_{\M(\R)}
    \end{align}
    is a Banach space. It is isometric to $\left(\mathcal{N}_{\Ls_1}\otimes\M(\R),\norm{\cdot}_{\bb{p}_1,\M(\R)}\right)$, under the mapping $\mathrm{proj}_{\mathcal{N}_{\Ls_1}}\otimes\Ls_2$ whose inverse is $\mathrm{proj}_{\mathcal{N}_{\Ls_1}}\otimes \Ls_{\bm{\phi}_2}^{-1}.$ 

    \item [3.]  The vector space $(iii)^{\star}$ expressed as $\left\{\sum_{n=1}^{N_2}\Ls_{\bm{\phi}_1}^{-1}\{m_{2,n}\}\otimes p_{2,n},\quad(m_{2,n})_{n=1}^{N_2}\in\M(\R)^{N_2}\right\}$ and equipped with the norm
    \begin{align}
        &\norm{\left[\Ls_1\otimes\mathrm{proj}_{\mathcal{N}_{\Ls_2}}\right]\left\{\sum_{n=1}^{N_2}\Ls_{\bm{\phi}_1}^{-1}\{m_{2,n}\}\otimes p_{2,n}\right\}}_{\M(\R),\bb{p}_2}\nonumber\\
        =&\norm{\sum_{n=1}^{N_2}m_{2,n}\otimes p_{2,n}}_{\M(\R), \bb{p}_2}:=\sum_{n=1}^{N_2}\norm{m_{2,n}}_{\M(\R)}
    \end{align}
    is a Banach space. It is isometric to $\left(\M(\R)\otimes\mathcal{N}_{\Ls_2},\norm{\cdot}_{\M(\R),\bb{p}_2}\right)$, under the mapping $\Ls_1\otimes\mathrm{proj}_{\mathcal{N}_{\Ls_2}}$ whose inverse is $\Ls_{\bm{\phi}_1}^{-1}\otimes\mathrm{proj}_{\mathcal{N}_{\Ls_2}}.$  

    \item  [4.] The vector space 
    $(iv)^{\star}$ expressed as 
    \begin{equation}
        \left\{\sum_{n,n'=1,1}^{N_1,N_2}d_{n,n'}p_{1,n}\otimes p_{2,n'},\quad (d_{n,n'})_{n,n'=1,1}^{N_1,N_2}\in\R^{N_1\times N_2}\right\}
    \end{equation}
     and equipped with the norm
\begin{align}
    &\norm{\sum_{n,n'=1,1}^{N_1,N_2}d_{n,n'}p_{1,n}\otimes p_{2,n'}}_{\bb{p}_1,\bb{p}_2}:=\sum_{n,n'=1,1}^{N_1,N_2}\vert d_{n,n'}\vert
\end{align}
is a Banach space.
\end{itemize}
\end{Proposition}

In Theorem \ref{th:native}, whose proof is in \ref{app:4.2}, we make a first step toward the description of the geometry of $\M_{\LL}(\R^2)$ by revealing that it is a Banach space for the norm $\norm{\cdot}_{\M_{\LL}}$. In addition, we show that $\M_{\LL}(\R^2)$ has a direct-sum structure, which will be central in the description of $[\LL]$-splines.

\begin{Theorem}
\label{th:native}The spaces $(i)^{\star},(ii)^{\star},(iii)^{\star}$, and $(iv)^{\star}$ are in direct sum and every $f\in\M_{\LL}(\R^2)$ has a unique direct-sum representation as 
\begin{align}
\label{eq:4.2.78}
f=&\left[\Ls_{\bm{\phi}_1}^{-1}\otimes \Ls_{\bm{\phi}_2}^{-1}\right]\{m\}+\sum_{n=1}^{N_1} p_{1,n}\otimes \Ls_{\bm{\phi}_2}^{-1}\{m_{1,n}\}+\sum_{n=1}^{N_2}\Ls_{\bm{\phi}_1}^{-1}\{m_{2,n}\}\otimes p_{2,n}\nonumber\\
&+\sum_{n,n'=1}^{N_1,N_2}d_{n,n'}p_{1,n}\otimes p_{2,n},
\end{align}
with $m\in \M(\R^2),$  $(m_{1,n})_{n=1}^{N_1}\in\M(\R)^{N_1},$ $(m_{2,n})_{n=1}^{N_2}\in\M(\R)^{N_2},$ and $(d_{n,n'})_{n,n'=1,1}^{N_1,N_2}\in\R^{N_1\times N_2}$. In addition, the mappings
\begin{equation}
    \Ls_1\otimes\Ls_2,\quad\mathrm{proj}_{\mathcal{N}_{\Ls_1}}\otimes \Ls_2,\quad\Ls_1\otimes\mathrm{proj}_{\mathcal{N}_{\Ls_2}},\quad \mathrm{proj}_{\mathcal{N}_{\Ls_1}}\otimes\mathrm{proj}_{\mathcal{N}_{\Ls_2}}§
\end{equation}
extend continuously from their domain of definition to the whole native space $\M_{\LL}(\R^2)$.
Consequently, $\M_{\LL}(\R^2)$ is a Banach space for the norm
\begin{align}
    \norm{f}_{\M_{\Ls_1\otimes\Ls_2}}&=\norm{[\Ls_1\otimes\Ls_2]\{f\}}_{\M(\R^2)}+\norm{\left[\mathrm{proj}_{\mathcal{N}_{\Ls_1}}\otimes \Ls_2\right]\{f\}}_{\bb{p}_1,\M(\R)}\nonumber\\
    &\quad+\norm{\left[\Ls_1\otimes\mathrm{proj}_{\mathcal{N}_{\Ls_2}}\right]\{f\}}_{\M(\R),\bb{p}_2}+\norm{\left[ \mathrm{proj}_{\mathcal{N}_{\Ls_1}}\otimes\mathrm{proj}_{\mathcal{N}_{\Ls_2}}\right]\{f\}}_{\bb{p}_1,\bb{p}_2}\nonumber\\
    &=\norm{m}_{\M(\R^2)}+\sum_{n=1}^{N_1}\norm{m_{1,n}}_{\M(\R)}+\sum_{n=1}^{N_2}\norm{m_{2,n}}_{\M(\R)}+\norm{\bb{d}}_{1}.\label{eq:4.2.79}
\end{align}
\end{Theorem}

We also introduce the seminorm \eqref{def:seminorm} as the sum of three norms, each one being effective on the spaces $(i)^{\star}$,$(ii)^{\star}$, and $(iii)^{\star}$, respectively. This seminorm is 
\begin{align} 
\snorm{\cdot}_{\M_{\Ls_1\otimes\Ls_2}}:=&\norm{[\Ls_1\otimes\Ls_2]\{\cdot\}}_{\M(\R^2)}+\norm{\left[\mathrm{proj}_{\mathcal{N}_{\Ls_1}}\otimes \Ls_2\right]\{\cdot\}}_{\bb{p}_1,\M(\R)}\nonumber\\
&+\norm{\left[\Ls_1\otimes\mathrm{proj}_{\mathcal{N}_{\Ls_2}}\right]\{\cdot\}}_{\M(\R),\bb{p}_2} \label{def:seminorm},
\end{align}
where $\norm{\cdot}_{\bb{p}_1,\M(\R)}$ and $\norm{\cdot}_{\M(\R),\bb{p}_2}$ are the mixed norms defined in Proposition \ref{prop:isom2}. We emphasize that
$\mathrm{proj}_{\mathcal{N}_{\Ls_i}}$ depends on the system $(\bb{p}_i,\bm{\phi}_i)$ and that, therefore, so does $\snorm{\cdot}_{\M_{\LL}}.$  We show in Proposition \ref{prop:4.4.1.1} that our native space is maximally defined for the seminorm $\snorm{\cdot}_{\M_{\LL}}$. The proof is given in \ref{app:4.2}.

\begin{Proposition}
\label{prop:4.4.1.1}
The vector space $\M_{\LL}(\R^2)$ admits the following definition
    \begin{equation}
    \label{eq:4.2.54}
        \M_{\LL}(\R^2)=\bigg\{f\in\D'(\R^2):\begin{cases}[\LL]\{f\}\in\M(\R^2)\\
        [\Ls_1\otimes\mathrm{proj}_{\mathcal{N}_{\Ls_2}}]\{f\}\in\M(\R)\otimes\mathcal{N}_{\Ls_2}\\
        [\mathrm{proj}_{\mathcal{N}_{\Ls_1}}\otimes\Ls_2]\{f\}\in \mathcal{N}_{\Ls_1}\otimes\M(\R)
        \end{cases}
        \bigg\},
    \end{equation}
which is equivalent to the one in Theorem \ref{th:native}.
\end{Proposition}
We observe that \eqref{eq:4.2.54} can be re-written in the more compact form
\begin{equation}
        \M_{\LL}(\R^2)=\left\{f\in\D'(\R^2):\quad\snorm{f}_{\M_{\LL}}<\infty\right\}.
\end{equation}
The seminorm $\snorm{\cdot}_{\M_{\LL}}$ not only defines the set $\M_{\LL}(\R^2)$, but also determines its geometry up to the choice of boundary conditions or, equivalently, the choice of admissible system. Concretely, the geometry is determined by the unit ball of $\snorm{\cdot}_{\M_{\LL}}$, itself characterized by its extreme points. We now specify these extreme points in Proposition \ref{prop:extremepoints}, whose proof is in \ref{app:4.2}.
\begin{Proposition}
\label{prop:extremepoints}
    The extreme points $e(\cdot)$ of the closed unit ball $\{f\in\M_{\LL}(\R^2)/\mathcal{N}_{\LL}:\quad\snorm{f}_{\M_{\LL}}\leq1\}$ of the norm $\snorm{\cdot}_{\M_{\LL}}$ are exactly represented by a function
    \begin{itemize}
        \item [1.]of the form $\pm\left[\Ls_{\bm{\phi}_1}^{-1}\otimes \Ls_{\bm{\phi}_2}^{-1}\right]\{\delta_{\bb{x}}\}=\pm g_{\bm{\phi}_1}(\cdot,x_1)\otimes g_{\bm{\phi}_2}(\cdot,x_2)$ for some $\bb{x}=(x_1,x_2)\in\R^2;$
        \item [2.]of the form $\pm p_{1,n}\otimes \Ls_{\bm{\phi}_2}^{-1}\{\delta_{x}\}=\pm p_{1,n}(\cdot)\otimes g_{\bm{\phi}_2}(\cdot,y)$ for some $y\in\R$, and $1\leq n\leq N_1;$
        \item [3.]of the form $\pm \Ls_{\bm{\phi}_1}^{-1}\{\delta_{x}\}\otimes p_{2,n}=\pm g_{\bm{\phi}_1}(\cdot,z)\otimes p_{2,n}(\cdot)$ for some $z\in\R$, and $1\leq n\leq N_2.$
    \end{itemize}    
\end{Proposition}
In the context of non-reflexive Banach spaces, the strong topology is too rigid for optimization purposes. In Theorem \ref{Th:7}, whose proof is given in \ref{app:4.2}, we establish the existence of the $\w$ topology.

\begin{table}[h!]
    \centering
    \begin{tabular}{c|c c}
    \hline\hline
       Properties $\backslash$ System Type  & $[\LL]$-admissible & $\bb{K}$-localized \\ \hline
         Strong Topology& system specific & system specific\\
         Geometry& system specific & system specific \\
         Solution-Set Size& suboptimal & optimal \\
         \hline\hline
    \end{tabular}
    \caption{Properties of the Banach space $\left(\M_{\LL}(\R^2),\norm{\cdot}_{\M_{\LL}}\right)$ in function of the type of $[\LL]$-admissible system. Based on the proof of Proposition \ref{prop:normequivalence} below, we hypothesise that the strong topology is system-specific. The geometry is described in Proposition \ref{prop:extremepoints}. The optimality of the size of the solution set is discussed in Section \ref{sec:5.3} and is based on the adequate weight distribution discussed in Section \ref{sec:fund}.}
    \label{tab:1}
\end{table}

\begin{Theorem}
\label{Th:7}The Banach spaces $(i)^{\star},(ii)^{\star},(iii)^{\star}$, and $(iv)^{\star}$ are the dual spaces of $(i),(ii),(iii)$, and $(iv),$ respectively. Consequently, $\M_{\Ls_1\otimes\Ls_2}(\R^2)$ equipped with the norm \eqref{eq:4.2.79} is the dual space of $\C_{\Ls_1\otimes\Ls_2}(\R^2)$ equipped with the norm \eqref{eq:4.1.44}.
\end{Theorem}

In the light of Theorem \ref{th:3.1} and \ref{Th:7}, it appears that $\M_{\LL}(\R^2)$ is the $\w$ completion of $\M_{\Ls_1}(\R)\otimes\M_{\Ls_2}(\R)$. Furthermore, the independence of $\C_{\LL}(\R^2)$, as a set, on the choice of admissible system does transfer by duality to $\M_{\LL}(\R^2)$. Consequently, the $\w$ topology does not depend on the choice of system.  We conclude this section with two generalizations.

\emph{Admissible Systems.} As in the one-dimensional case, the characterization of the predual $\C_{\LL}(\R^2)$ enables us to extend the notion of $[\LL]$-admissible system to a collection  $((\phi_{1,n})_{n=1}^{N_1},(\phi_{2,n'})_{n'=1}^{N_2})$ with $\phi_{1,n}\in\C_{\Ls_1}(\R)$ and $\phi_{2,n}\in\C_{\Ls_2}(\R)$. In particular, one has that
\begin{equation}
    \phi_{1,n}\otimes\phi_{2,n}\in\C_{\Ls_1}(\R)\otimes\C_{\Ls_2}(\R)\subset\C_{\LL}(\R^2).
\end{equation}
Indeed, we observe that the constructions and results in Sections \ref{sec:4.1} and \ref{sec:4.2} hold if one replaces, a posteriori, a universal admissible system by an admissible one. In the sequel, we shall implicitly assume that $((\bb{p}_1, \bb{p}_2), (\bm{\phi}_1,\bm{\phi}_2))$ is only an $[\LL]$-admissible system.

\emph{K-localized Systems.} General $[\LL]$-admissible systems give rise to the same inconvenient integral representation as in the one-dimensional case. As in Section \ref{sec:fund}, a first remedy is to define localized systems. Let $\bb{K}=(\K_1\times\K_2)$ with $\K_1$ and $\K_2$ two intervals. We say that an $[\LL]$-admissible system $((\bb{p}_1, \bb{p}_2), (\bm{\phi}_1,\bm{\phi}_2))$ is $\bb{K}$-localized if $(\bb{p}_1,\bm{\phi}_1)$ and $(\bb{p}_2,\bm{\phi}_2)$ are $\K_1$-localized and $\K_2$-localized. This topic, along with its extension to fundamental systems, will be discussed in greater depth in Section \ref{sec:4.5}.

Finally, we summarize some key aspects of $\left(\M_{\LL}(\R^2),\norm{\cdot}_{\M_{\LL}}\right)$ in Table \ref{tab:1}, that will be central in the resolution of the optimization problem in Section \ref{sec:5}.

\subsection{Splines, Regularity, and Growth}
\label{sec:4.3}
\begin{Definition}
\label{def:spline}
    A distribution $f\in\D'(\R^2)$ is an $[\LL]$-spline if,  $\forall\bb{t}=(t_1,t_2)\in\R^2$, it can be written as
    \begin{align}
    f(t_1,t_2)=&\sum_{k=1}^Ka_{k}g_{\Ls_1}(t_1-x_{1,k}) g_{\Ls_2}(t_2-x_{2,k})+\sum_{n=1}^{N_1}\sum_{m=1}^{K_{1,n}}b_{n,m}p_{1,n}(t_1) g_{\Ls_2}(t_2-y_{n,m})\nonumber\\
    &+\sum_{n'=1}^{N_2}\sum_{m'=1}^{K_{2,n'}}c_{n',m'}g_{\Ls_1}(t_1-z_{n',m'}) p_{2,n'}(t_2)+\sum_{n,n'=1,1}^{N_1,N_2}d_{n,n'}p_{1,n}(t_1) p_{2,n'}(t_2),
    \end{align}    
    where $a_k\in\R,b_{n,m}\in\R,$ and $c_{n',m'}\in\R$ are the weights and $\bb{x}_k=(x_{1,k},x_{2,k})\in\R^2,y_{n,m}\in\R,$ and $z_{n',m'}\in\R$ are the knots with $\left(K+\sum_{n=1}^{N_1}K_{1,n}+\sum_{n'=1}^{N_2}K_{2,n'}\right)$ being the number of knots.
\end{Definition}
If $f_1$ is an $\Ls_1$-spline with $K_1$ knots and $f_2$ is an $\Ls_2$-spline with $K_2$ knots, then $f_1\otimes f_2$ is an $[\LL]$-spline with at most $( K_1K_2+N_1K_2+N_2K_1)$ knots. Thus, the tensor product of two one-dimensional splines is a particular case of a tensor-product-spline, with a  specific knot structure. We provide in Proposition \ref{prop:4.4.1.2} an equivalent definition of $[\LL]$-splines based on the technique of Proposition \ref{prop:4.4.1.1}, from which the proof follows directly.
\begin{Proposition}
\label{prop:4.4.1.2}
    A distribution $f\in\D'(\R^2)$ is an $[\LL]$-spline if and only if it holds jointly that
    \begin{align}
        \left[\Ls_1\otimes\Ls_2\right]\{f\}&=\sum_{k=1}^Ka_k\delta_{\bb{x}_k}\nonumber\\
        [\mathrm{proj}_{\mathcal{N}_{\Ls_1}}\otimes\Ls_2]\{f\}&=\sum_{n=1}^{N_1}p_{1,n}\otimes\left( \sum_{m=1}^{K_{1,n}}b_{n,m}\delta_{y_{n,m}}\right)\nonumber\\
        [\Ls_1\otimes\mathrm{proj}_{\mathcal{N}_{\Ls_2}}]\{f\}&=\sum_{n'=1}^{N_2}\left( \sum_{m'=1}^{K_{2,n'}}c_{n',m'}\delta_{z_{n',m'}}\right)\otimes p_{2,n'}\label{eq:4.3.51}.
    \end{align}  
\end{Proposition}

Accordingly, we can interpret a spline as a finite linear combination of extreme points of the unit ball in $\M_{\LL}(\R^2)$, themselves specified by the differential operator $\LL$. In Propositions \ref{prop:regularity}, we show that the regularity properties, inherited from the order of $\LL$, transfer to the whole space $\M_{\LL}(\R).$ The proof of which is given in \ref{app:4.3}.
\begin{Proposition}
\label{prop:regularity}
Let $N_1,N_2$ be the degrees of the ODO $\Ls_1,\Ls_2$. Then, $\forall f\in\M_{\LL}(\R^2),$
    \begin{align}
        &\forall d_1\in[0\ldots(N_1-1)],\forall d_2\in[0\ldots(N_2-1)]:\quad\quad\frac{\partial^{d_1+d_2}f}{\partial t_1^{d_1}\partial t_2^{d_2}}\in\mathcal{L}_{\infty,\mathrm{loc}}(\R^2)\\
        &\forall d_1\in[0\ldots(N_1-2)],\forall d_2\in[0\ldots(N_2-2)]:\quad\quad\frac{\partial^{d_1+d_2}f}{\partial t_1^{d_1}\partial t_2^{d_2}}\in\C(\R^2).
    \end{align}
\end{Proposition}

Likewise, it is straightforward (but lengthy) to prove that the growth of the function $f\in\M_{\LL}(\R^2)$ is of order $\mathcal{O}(t_i^{N_i-1}\e^{\alpha_it})$ in the dimension $i\in\{1,2\}$.

\subsection{Search-Space Localization}
\label{sec:4.5}
Following the same path as in the one-dimensional case, we set  
\begin{equation}
    \bb{K}=\mathrm{K}_1\times\mathrm{K}_2=[K_1^-,K_1^+]\times[K_2^-,K_2^+]
\end{equation}
and define the Banach subspace $\M_{\Ls_1\otimes\Ls_2}(\bb{K})\subset\M_{\LL}(\R^2)$ as the localized search space
\begin{align}
    \M_{\Ls_1\otimes\Ls_2}(\bb{K})=&\underbrace{\left([\Ls_{\bm{\phi}_1}^{-1}\otimes \Ls_{\bm{\phi}_2}^{-1}]\{\M(\bb{K})\}\right)}_{(i)^{\star}}\oplus\underbrace{\left(\mathcal{N}_{\Ls_1}\otimes \Ls_{\bm{\phi}_2}^{-1}\{\M(\mathrm{K}_2)\}\right)}_{(ii)^{\star}}\nonumber\\
&\oplus\underbrace{\left(\Ls_{\bm{\phi}_1}^{-1}\{\M(\mathrm{K}_1)\}\otimes\mathcal{N}_{\Ls_2}\right)}_{(iii)^{\star}}\oplus\underbrace{\left(\mathcal{N}_{\Ls_1}\otimes\mathcal{N}_{\Ls_2}\right)}_{(iv)^{\star}}.\label{eq:4.5.55}
\end{align}
Observe that $\M_{\LL}(\bb{K})$ is the $\w$ completion of $\M_{\Ls_1}(\mathrm{K}_1)\otimes\M_{\Ls_2}(\mathrm{K}_2)$ and that it admits the equivalent definition 
\begin{equation}
    \M_{\LL}(\bb{K})=\bigg\{f\in\D'(\R^2):\begin{cases}[\LL]\{f\}\in\M(\bb{K})\\
    [\Ls_1\otimes\mathrm{proj}_{\mathcal{N}_{\Ls_2}}]\{f\}\in\M(\K_1)\otimes\mathcal{N}_{\Ls_2}\\
    [\mathrm{proj}_{\mathcal{N}_{\Ls_1}}\otimes\Ls_2]\{f\}\in \mathcal{N}_{\Ls_1}\otimes\M(\K_2)
    \end{cases}\bigg\}.
\end{equation}

We show in Proposition \ref{prop:CompactRepresentation}, whose proof is in \ref{app:4.5}, that in contrast to functions in $\M_{\LL}(\R^2)$, those in $\M_{\LL}(\bb{K})$ enjoy a simple integral representation.
\begin{Proposition}
    \label{prop:CompactRepresentation}
    A distribution $f\in\D'(\R^2)$ belongs to the space $\M_{\LL}(\bb{K})$ if and only if it is of the form, $\forall \bb{t}=(t_1,t_2)\in\R^2$,
    \begin{align}
        f(\bb{t})=&(g_{\Ls_1}\otimes g_{\Ls_2}\ast m)(\bb{t})+\sum_{n=1}^{N_1}p_{1,n}(t_1)(g_{\Ls_2}\ast m_{1,n})(t_2)\nonumber\\
        &+\sum_{n=1}^{N_2}(g_{\Ls_1}\ast m_{2,n})(t_1)p_{2,n}(t_2)
        +\sum_{n,n'=1,1}^{N_1,N_2}d_{n,n'}p_{1,n}(t_1)p_{2,n'}(t_2) \label{eq:4.5.51},
\end{align}
with $m\in\M(\bb{K}),(m_{1,n})_{n=1}^{N_1}\in\M(\mathrm{K}_2)^{N_1},(m_{2,n})_{n=1}^{N_2}\in\M(\mathrm{K}_1)^{N_2}$, and $(d_{n,n'})_{n,n'=1,1}^{N_1,N_2}\in\R^{N_1\times N_2}.$
\end{Proposition}

Equation \eqref{eq:4.5.51} can be refined into a direct-sum representation if one considers the appropriate  type of $[\LL]$-admissible system. Let $(\bb{p}_i,\bm{\iota}_i)$ be the unique $\K_i$-fundamental system. We say that the concatenation $((\bb{p}_1,\bb{p}_2),(\bm{\iota}_1,\bm{\iota}_2))$ of two $\bb{K}_i$-fundamental systems is a fundamental $\bb{K}$-localized system. In virtue of Proposition, \ref{prop:ideal2} it is unique.

This choice of fundamental system leads to a simplification of the kernels, from $g_{\bm{\iota}_i}$ to $g_{\Ls_i}$, and from $g_{\bm{\iota}_1}\otimes g_{\bm{\iota}_2}$ to $g_{\Ls_1}\otimes g_{\Ls_2}$. In addition, one can use $\Ls^{-1}_i$ instead of $\Ls^{-1}_{\bm{\iota}_i}$, and $\Ls^{-1}_1\otimes\Ls^{-1}_2$ instead of $\Ls^{-1}_{\bm{\iota}_1}\otimes\Ls^{-1}_{\bm{\iota}_2}$, where the integral operator $\Ls^{-1}_1\otimes\Ls_2^{-1}$ is described by the kernel $g_{\Ls_1}\otimes g_{\Ls_2}$ and 
\begin{equation}
    \left[\Ls_1^{-1}\otimes\Ls_2^{-1}\right]\{\M(\bb{K})\}=\left\{\int_{\bb{K}}g_{\Ls_1}(\cdot-x_1)g_{\Ls_2}(\cdot-x_2)\mathrm{d}m(x_1,x_2): m\in\M(\bb{K})\right\}\label{eq:4.4.75}.
\end{equation}

These simplified representations can be injected in \eqref{eq:4.5.55} to establish that
\begin{align}
    \M_{\Ls_1\otimes\Ls_2}(\bb{K})=&\underbrace{\left([\Ls_{1}^{-1}\otimes \Ls_{2}^{-1}]\{\M(\bb{K})\}\right)}_{(i)^{\star}}\oplus\underbrace{\left(\mathcal{N}_{\Ls_1}\otimes \Ls_{2}^{-1}\{\M(\mathrm{K}_2)\}\right)}_{(ii)^{\star}}\nonumber\nonumber\\
&\oplus\underbrace{\left(\Ls_{1}^{-1}\{\M(\mathrm{K}_1)\}\otimes\mathcal{N}_{\Ls_2}\right)}_{(iii)^{\star}}\oplus\underbrace{\left(\mathcal{N}_{\Ls_1}\otimes\mathcal{N}_{\Ls_2}\right)}_{(iv)^{\star}}.
\end{align}

%\begin{wrapfigure}{r}{0.55\textwidth}
%\vspace{-0.5cm}
 % \begin{center}
  %  \includegraphics[width=0.48\textwidth]{distribution.png}
 % \end{center}
 % \caption{Schematization of the spatial influence of functions in $(i)^{\star}$, $(ii)^{\star}$ and $(iii)^{\star}$.}
%\end{wrapfigure}

\begin{figure}
    \centering
    \includegraphics[width=0.5\linewidth]{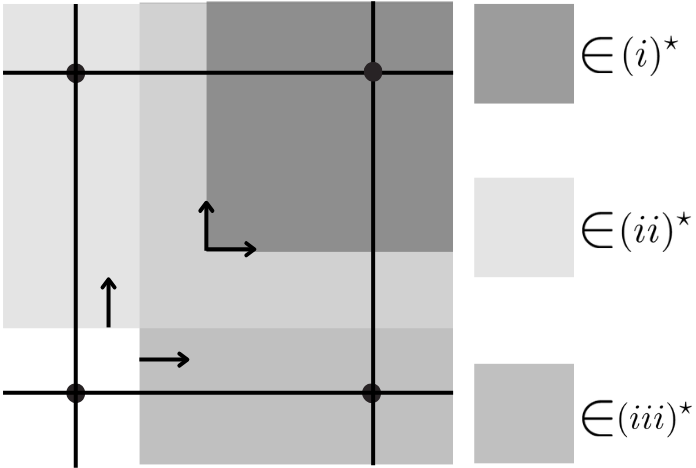}
  \caption{Schematization of the spatial influence of functions in $(i)^{\star}$, $(ii)^{\star}$, and $(iii)^{\star}$.}
\end{figure}
We note that, the extreme points in Proposition \ref{prop:extremepoints} also enjoy this kernel simplification. After this geometrical update, we study in Proposition \ref{prop:normequivalence}, whose proof is in \ref{app:4.5}, the seemingly different strong topologies on  $\M_{\LL}(\bb{K})$.

\begin{Proposition}
\label{prop:normequivalence}
     On $\M_{\LL}(\bb{K})$, the norm generated by an arbitrary $[\LL]$-admissible system and the norm generated by the $\bb{K}$-fundamental system are equivalent.
\end{Proposition}

Proposition \eqref{prop:normequivalence} reveals that the strong topology does not depend on the choice of system. Therefore, without loss of generality, one can use use the favourable $\bb{K}$-fundamental system. Relevant properties of $\left(\M_{\LL}(\bb{K}),\norm{\cdot}_{\M_{\LL}}\right)$ are summarized in Table \ref{tab:2}. The equivalence of the strong topologies is established in Proposition \ref{prop:normequivalence}. The geometry description follows from Proposition \ref{prop:extremepoints}. The optimality of the size of the solution set is discussed in Section \ref{sec:5.3}.

\begin{table}[h!]
    \centering
    \begin{tabular}{c|c c c}
    \hline\hline
       Properties $\backslash$ System Type  & $[\LL]$-admissible & $\bb{K}$-localized & $\bb{K}$-Fundamental \\ \hline
         Strong Topology& equivalent & equivalent & equivalent\\
         Geometry& system-specific & system-specific & \makecell{equivalent}\\
         Solution-Set Size& suboptimal & optimal & optimal\\
         \hline\hline
    \end{tabular}
    \caption{Properties of the Banach space $\left(\M_{\LL}(\bb{K}),\norm{\cdot}_{\M_{\LL}}\right)$ in function of the type of $[\LL]$-admissible system.}
    \label{tab:2}
\end{table}

%%%%%%%%%%%%%%%%%%%%%%%%%%%%%%%%%%%%%%%%%%%%%%%%%%%%%
\section{Main Results: Resolution of Two-Dimensional Continuous-Domain Inverse Problems}
\label{sec:5}

\subsection{Assumptions}
\label{sec:5.1}
\begin{itemize}
    \item []\hspace{-0.5cm}\emph{Assumption 1} The operators $\Ls_1=(\Dd-\alpha_1\mathrm{I})^{N_1}$ and $\Ls_2=(\Dd-\alpha_2\mathrm{I})^{N_2}$ are ODO.
    \item []\hspace{-0.5cm}\emph{Assumption 2} \label{ass:1}The measurement operator $\langle\cdot,\bm{\nu}\rangle:\M_{\Ls_1\otimes\Ls_2}(\R^2)\to\R^M$ is linear, $\mathrm{weak}^{\star}$-continuous, surjective, and injective on $\mathcal{N}_{\Ls_1}\otimes\mathcal{N}_{\Ls_2}$.
    \item []\hspace{-0.5cm}\emph{Assumption 3} \label{ass:4} The data-fidelity functional $E:\R^M\times\R^M\to\R^+\cup\{\infty\}$ is proper, coercive, lower semicontinuous, and strictly convex in its second argument.
    \item []\hspace{-0.5cm}\emph{Assumption 4} The system $((\bb{p}_1,\bb{p}_2),(\bm{\phi}_1,\bm{\phi}_2))$ is $[\LL]$-admissible.
    \item []\hspace{-0.5cm}\emph{Assumption 5} The measurement operator $\bm{\nu}$ is compactly supported in $\bb{K}=(\mathrm{K}_1\times\mathrm{K}_2)=\big([K_1^-,K_1^+]\times[K_2^-,K_2^+]\big)$.

    \item []\hspace{-0.5cm}\emph{Assumption 6} The system $((\bb{p}_1,\bb{p}_2),(\bm{\iota}_1,\bm{\iota}_2))$ is the unique ideal $\bb{K}$-localized system, given in Proposition \ref{prop:ideal2}.
\end{itemize}

\subsection{Representer Theorem for General Measurement Operators}
\label{sec:5.2}

In Theorem \ref{th:optweakstar} we reveal that for the search space $\mathcal{X}=\M_{\LL}(\R^2)$, the extreme points $f^{\star}$ of the solution  set $\mathcal{V}$ are the linear combinations of kernels $g_{\bm{\phi}_1},g_{\bm{\phi}_2}$ and $g_{\bm{\phi}_1}\otimes g_{\bm{\phi}_2}$ with at most $(M-N_1N_2)$ knots.

\begin{Theorem}
\label{th:optweakstar}
If Assumptions 1 to 4 are verified, then, for $\mathcal{X}=\M_{\LL}(\R^2)$, the solution set $\mathcal{V}$ is nonempty, $\w$-compact, and is the $\w$-closed convex hull of its extreme points, which all are of the form
\begin{align}
\label{eq:5.2.66}
    f^{\star}(t_1,t_2)=&\sum_{k=1}^Ka_{k}g_{\bm{\phi}_1}(t_1,x_{1,k}) g_{\bm{\phi_2}}(t_2,x_{2,k})\nonumber\\
    &+\sum_{n=1}^{N_1}\sum_{m=1}^{K_{1,n}}b_{n,m}p_{1,n}(t_1) g_{\bm{\phi}_2}(t_2,y_{n,m})\nonumber\\
    &+\sum_{n'=1}^{N_2}\sum_{m'=1}^{K_{2,n'}}c_{n',m'}g_{\bm{\phi}_1}(t_1,z_{n',m'}) p_{2,n'}(t_2)\nonumber\\
    &+\sum_{n,n'=1,1}^{N_1,N_2}d_{n,n'}p_{1,n}(t_1) p_{2,n'}(t_2),
\end{align}
with 
\begin{equation}
K+\sum_{n=1}^{N_1}K_{1,n}+\sum_{n'=1}^{N_2}K_{2,n'}\leq M-N_1N_2,
\end{equation}
and $(x_{1,k},x_{2,k})\in\R^{2},y_{n,m}\in\R$ and $z_{n',m'}\in\R.$
\end{Theorem}

\begin{proof}[\textbf{Proof of Theorem \ref{th:optweakstar}.}]
In virtue of Assumptions 1 to 4, we apply \citep[Proposition 8]{gupta2018continuous} to find that $\mathcal{V}$ is nonempty, convex, and $\w$-compact. Then, the Krein-Milman theorem yields that $\mathcal{V}$ is the $\w$-closed convex hull of its extreme points. Next, we apply \citep[Corollary 3.8]{boyer2019representer} to conclude that these extreme points all are linear combinations of at most $(M-N_1N_2)$ extreme points of the unit ball in $\left(\M_{\LL}(\R^2)/\mathcal{N}_{\LL},\snorm{\cdot}_{\M_{\LL}}\right)$. We conclude the proof with the extreme-point characterization of Proposition \ref{prop:extremepoints}.
\end{proof}

Theorem \ref{th:optweakstar} establishes the optimality of $[\LL]$-splines. There, the extreme-point solutions are described as the sum of four terms. 
\begin{itemize}
    \item [1.] The first term, formed as a tensor product of kernels, describes the component of the solution that is adaptive in both dimensions. This term is regularized, and its weights $a_k$ are controlled by the seminorm $\norm{[\LL]\{\cdot\}}_{\M(\R^2)}$ with
    \begin{equation}
        \norm{\left[\LL\right]\left\{\sum_{k=1}^Ka_{k}g_{\bm{\phi}_1}(\cdot,x_{1,k}) g_{\bm{\phi_2}}(\cdot,x_{2,k})\right\}}_{\M(\R^2)}=\sum_{k=1}^K\vert a_k\vert.
    \end{equation}
    \item [2\&3.] The second and third terms are formed as the tensor product between exponential polynomials and a kernel. For the second  term, only the component of the second dimension is adaptive, and its weights $b_{n,m}$ are controlled by 
        \begin{equation}
        \norm{\left[\mathrm{proj}_{\mathcal{N}_{\Ls_1}}\otimes\Ls_2\right]\left\{\sum_{n=1}^{N_1}\sum_{m=1}^{K_{1,n}}b_{n,m}p_{1,n}(\cdot) g_{\bm{\phi}_2}(\cdot,y_{n,m})\right\}}_{\bb{p}_1,\M(\R)}=\sum_{n=1}^{N_1}\sum_{m=1}^{K_{1,n}}\vert b_{n,m}\vert.
    \end{equation}
    \item [4.] The fourth term is the tensor product of exponential polynomials and is regularization-free. This component is viewed as a classic linear-regression term.
\end{itemize}
Nevertheless, Theorem \ref{th:optweakstar} is not fully satisfactory as the extreme points are expressed in the kernel basis which lacks the key property of shift-invariance. The lack of localization of the knots $(x_{1,k},x_{2,k})$, $y_{n,m}$, and $z_{n',m'}$ is also a concerning limitation. These two issues are resolved when one assumes $\bm{\nu}$ to be compactly supported.

\subsection{Representer Theorems for Compactly Supported Measurement Operators}
\label{sec:5.3}
Denote by $\overset{\circ}{\bb{K}},\overset{\circ}{\mathrm{K}}_1,\overset{\circ}{\mathrm{K}}_2$ the interiors of $\bb{K}$, $\mathrm{K}_1,\mathrm{K}_2,$ respectively. Theorem \ref{th:RTlocalized}, whose proof is in \ref{app:5.3}, gives the localization of the knots in \eqref{eq:5.2.66} when $\bm{\nu}$ is supported in $\bb{K}$.

\begin{Theorem}
\label{th:RTlocalized}
 If Assumptions 1 to 5 are verified, then the knots of the extreme points are localized, in the sense that $z_{n',m'}\in\mathrm{K}_1,y_{n,m}\in\mathrm{K}_2$ and $(x_{1,k},x_{2,k})\in\tilde{\bb{K}}$, with 
    \begin{equation}
        \tilde{\bb{K}}=[\mathrm{min}(\phi^-_1,K^{-}_1),\mathrm{max}(\phi^+_1,K^+_1)]\times[\mathrm{min}(\phi^-_2,K^{-}_2),\mathrm{max}(\phi^+_2,K^+_2)].
    \end{equation}
If in addition $((\bb{p}_1,\bb{p}_2),(\bm{\phi}_1,\bm{\phi}_2))$ is $\bb{K}$-localized, then $(x_{1,k},x_{2,k})\in\bb{K}$ and the solution set $\mathcal{V}$ is a subset of $\M_{\LL}(\bb{K})$ and, consequently,
     \begin{equation}
     \underset{f\in\M_{\Ls_1\otimes\Ls_2}(\R^2)}{\mathrm{argmin}}\quad \mathcal{J}(f)=\underset{f\in\M_{\Ls_1\otimes\Ls_2}(\bb{K})}{\mathrm{argmin}}\mathcal{J}(f).
 \end{equation}
\end{Theorem}

\begin{figure}
    \centering
    \includegraphics[width=0.5\linewidth]{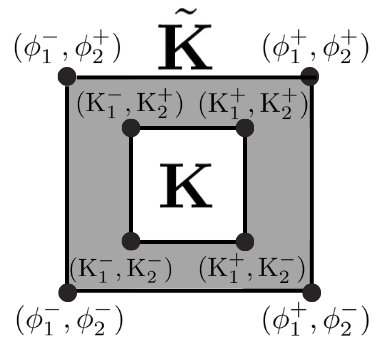}
\caption{ Schematic of the squares $\bb{K}$ and $\tilde{\bb{K}}$. The grey zone represents the additional area where a knot $(x_{1,k},x_{2,k})$ may be placed due to a suboptimal choice of admissible system.}
\label{fig:support}
\end{figure}

Theorem \ref{th:RTlocalized} shows that Assumption 5 is sufficient to localize the knots in a compact set $\tilde{\bb{K}}$, as illustrated in Figure \ref{fig:support}. In addition, this localization is optimal, in the sense that $\tilde{\bb{K}}=\bb{K}$ if and only if the admissible system is $\bb{K}$-localized. In this case, the optimisation of the functional $\mathcal{J}$ on $\M_{\LL}(\R^2)$ is equivalent to the optimisation on the simpler, compact-domain space $\M_{\LL}(\bb{K})$.

On the practical hand, the numerical resolution of the problem only needs to focus on the finding of knots inside $\bb{K}$, which is time-saving for the practitioner. On the theoretical hand we specialize Theorem \ref{th:RTlocalized} in Theorem \ref{th:RTcompact}, with the knowledge that, when the system is $\bb{K}$-localized, it is sufficient to optimize on $\M_{\LL}(\bb{K})$ to find all solutions. There, we show that the  $\bb{K}$-fundamental system yields extreme point solutions that are expressed in the Green's function basis. The proof is similar to the one of Theorem \ref{th:optweakstar} and is omitted.

\begin{Theorem}
\label{th:RTcompact}
If Assumptions 1 to 6 are verified and $\mathcal{X}=\M_{\LL}(\bb{K})$, then the solution set $\mathcal{V}$ is nonempty, $\mathrm{weak}^{\star}$-compact, and is the $\mathrm{weak}^{\star}$-closed convex hull of its extreme points which are $[\LL]$-splines \eqref{def:spline} with parameters that verify the following properties.
\begin{itemize}
    \item []\emph{Sparsity:} there are at most $(M-N_1N_2)$ knots. 
    \item []\emph{Localization:} the knots are localized in $\bb{K}$.
\end{itemize}
\end{Theorem}

In Corollary \ref{coro:interior}, we leverage the specific choice of basis from Item 3 of the $\K_i$-fundamental assumption to show that, without loss of generality, one can always assume that the knots of an extreme-point solution are localized in $\overset{\circ}{\bb{K}}$.

\begin{Corollary} 
\label{coro:interior}
If the assumptions of Theorem \ref{th:RTcompact} are verified and $f^{\star}$ is an extreme-point solution, then there exists another extreme-point solution $f^{\star\star}$ such that 
\begin{itemize}
    \item [1.] the knots $(x_{1,k},x_{2,k}),y_{n,m}$, and $z_{n',m'}$ are in $\overset{\circ}{\bb{K}},\overset{\circ}{\mathrm{K}}_1$, and $\overset{\circ}{\mathrm{K}}_2$, respectively;
    \item [2.] over the square $\bb{K}$, it holds that $f^{\star}(\bb{t})=f^{\star\star}(\bb{t})$;
    \item [3.] over the regularization, it holds that $\snorm{f^{\star}}_{\M_{\LL}}=\snorm{f^{\star\star}}_{\M_{\LL}}$.
\end{itemize} 
\end{Corollary}

\subsection{Role of Admissible Systems and Regularization on the Boundary}
\label{sec:5.4}

On the one hand, Theorem \ref{th:optweakstar} shows that, when $\bm{\nu}$ is arbitrary, \emph{sparsity} is preserved but \emph{localization} and \emph{shift-invariance} of the extreme-point solutions have to be sacrificed. In this case, there is no reason to prefer a specific $[\LL]$-admissible system. On the other hand, when $\bm{\nu}$ is supported in $\bb{K}$, we showed that the $\bb{K}$-fundamental system is the unique system that yields the sparsity, localization in $\overset{\circ}{\bb{K}}$, and shift invariance of the extreme-point solutions. Concretely, there exists a unique ``good" way of regularizing on the null space of $\LL$. It is to benefit from the seminorm
\begin{align}
\label{eq:5.4.77}
    \snorm{f}_{\M_{\LL}}=&\norm{[\LL]\{f\}}_{\M(\bb{K})}+\sum_{n=1}^{N_1}\norm{[(\Dd-\alpha_1\mathrm{I})^{n-1}\otimes\Ls_2]\{f\}}_{\M(\{K_{1}^-\}\times\mathrm{K}_2)}\nonumber\\
    &+\sum_{n'=1}^{N_2}\norm{[\Ls_1\otimes(\Dd-\alpha_2\mathrm{I})^{n'-1}]\{f\}}_{\M(\mathrm{K}_1\times\{K_2^-\})},
\end{align}
where we used the convention that the limit of a discontinuous function is by default the left-limit. The null-space regularization part can be seen as partial-derivative sliding on the bottom-left corner of $\partial\bb{K}.$

\subsection{Relation Between Regularization and the Principle of Causality}
We fix $\bb{K}=[0,1]^2$ and $\Ls_1=\Ls_1=\Dd^N$. Anchored in the one-dimensional theory, and the field of signal processing, we transposed in $D=2$ the principle of causality \citep{reichenbach1991direction}. Concretely, we chose there to work with the causal Green's function $g_{\Dd^N}\otimes g_{\Dd^N}$, and the causal $[0,1]^2$-fundamental system. In consequence,
\begin{itemize}
    \item the extreme-point solutions are described in terms of the causal Green's function;
    \item in the null-space regularization, the partial derivatives are only integrated on the lower part  $\big(\{0\}\times[0,1]\big)\cup\big([0,1]\times\{0\}\big)$ of the rectangle $\bb{K}$.
\end{itemize}
The same construction with the unique anti-causal Green's function  can be made. In this case, the partial derivatives in the regularization are integrated only on the upper-right corner $\big(\{1\}\times[0,1]\big)\cup\big([0,1]\times\{1\}\big)$. In order to regularize on the whole boundary $\partial[0,1]^2$, one may use the mixture
\begin{align}
    \snorm{f}_{\M_{\Dd^N\otimes\Dd^N}}=&\norm{[\Dd^N\otimes\Dd^N]\{f\}}_{\M([0,1]^2)}+\sum_{n=1}^{N}\norm{[\Dd^{n-1}\otimes\Dd^N]\{f\}}_{\M(\partial[0,1]\times[0,1])}\nonumber\\
    &+\sum_{n=1}^{N}\norm{[\Dd^N\otimes\Dd^{n-1}]\{f\}}_{\M([0,1]\times\partial[0,1])}\label{eq:5.5.100}.
\end{align}
Further work needs to be done to calculate the extreme-point solutions associated to this $\emph{acausal}$ regularization. Remark that, for $N=1$, \eqref{eq:5.5.100} simplifies into 
\begin{align}
    \snorm{f}_{\M_{\Dd\otimes\Dd}}=\norm{[\Dd\otimes\Dd]\{f\}}_{\M([0,1]^2)}&+\norm{\nabla_{\bb{p}}\{f\}}_{\M(\partial[0,1]^2)},
\end{align}
where $\bb{p}$ is the unit vector normal to the boundary.

%%%%%%%%%%%%%%%%%%%%%%%%%%%%%%%%%%%%%%%%%%%%%%%%%%%%%
\section{Examples}
\label{sec:6}

\subsection{Piecewise-Constant Splines}
We fix $\bb{K}=[0,1]^2$ and specify $\Ls_1=\Ls_2=\Dd$. The causal Green's function of $\Dd$ is the Heaviside function $(\cdot)_+^0$ and $\mathcal{N}_{\Dd}=\R$. Then, the localized native space is
\begin{align}
        \M_{\DD}([0,1]^2)=&\left([\Dd^{-1}\otimes\Dd^{-1}]\{\M([0,1]^2)\}\right)\oplus\left(\R\otimes \Dd^{-1}\{\M([0,1])\}\right)\nonumber\\
        &\oplus\left(\Dd^{-1}\{\M([0,1])\}\otimes\R\right)\oplus\R,
\end{align}
and $f\in\M_{\DD}([0,1]^2)$ admits the direct-sum decomposition
\begin{align}
    f(\bb{t})=&\int_{[0,1]^2}(t_1-x_1)_+^0(t_2-x_2)_+^0\mathrm{d}m(x_1,x_2)&+\int_{[0,1]}(t_1-x_1)_+^0\mathrm{d}m(x_1)\nonumber\\
    &+\int_{[0,1]}(t_2-x_2)_+^0\mathrm{d}m(x_2)+d_1
\end{align}
 with $m\in\M([0,1]^2),(m_1,m_2)\in\M([0,1])^2$, and $d_1\in\R$. Let $\snorm{\cdot}_{\M_{\DD}}$ be the seminorm defined in \eqref{eq:5.4.77}. If Assumptions 1 to 6 are verified, then 
\begin{equation}
     \mathcal{V}=\underset{f\in\M_{\DD}([0,1]^2)}{\text{argmin}}\Big( E(\bb{y},\langle\ f,\bm{\nu}\rangle)+\lambda\snorm{f}_{\M_{\DD}}\Big)
\end{equation}
is nonempty, $\mathrm{weak}^{\star}$-compact, and is the $\mathrm{weak}^{\star}$-closed convex of its extreme points which are of the form 
\begin{align}
\label{eq:6.1.68}
    \hspace{-0.3cm}f(\bb{t})=&\sum_{k=1}^Ka_{k}(t_1-x_{1,k})_+^0 (t_2-x_{2,k})_+^0+\sum_{m=1}^{K_{1}}b_{m}(t_1-y_{m})_+^0\nonumber\\&+\sum_{m'=1}^{K_{2}}c_{m'}(t_2-z_{m'})_+^0 +d_1,
\end{align}
with the number of knots $\left(K+K_1+K_2\right)$ being upper-bounded by $\left( M-1\right)$. Functions of the form \ref{eq:6.1.68} are illustrated in Figures \ref{fig:1} and \ref{fig:2}.

\subsection{Piecewise-Linear Splines}
We specify $\Ls_1=\Ls_2=\Dd^2$. Recall that the causal Green's function of $\Dd^2$ is the ReLU function $(\cdot)_+$ and that $\mathcal{N}_{\Dd^2}$ is the vector space of linear functions. Then, the localized native space is 
\begin{align}
        \M_{\DDD}([0,1]^2)=&\left([\Dd^{-2}\otimes\Dd^{-2}]\{\M([0,1]^2)\}\right)\oplus\left(\mathcal{N}_{\Dd^2}\otimes\Dd^{-2}\{\M([0,1])\}\right)\nonumber\\
        &\oplus\left(\Dd^{-2}\{\M([0,1])\}\otimes\mathcal{N}_{\Dd^2}\right)\oplus\left(\mathcal{N}_{\Dd^2}\otimes\mathcal{N}_{\Dd^2}\right),
\end{align}
where $f\in\M_{\DDD}([0,1]^2)$ admits the direct-sum decomposition 
\begin{align}
    f(\bb{t})=&\int_{[0,1]^2}(t_1-x_1)_+(t_2-x_2)_+\mathrm{d}m(x_1,x_2)+\sum_{n=1}^{2}\left(\int_{[0,1]}(t_2-x_2)_+\mathrm{d}m_{1,n}(x_2)\right)t_1^{n-1}\nonumber\\
    &+\sum_{n=1}^{2}\left(\int_{[0,1]}(t_1-x_1)_+\mathrm{d}m_{2,n}(x_1)\right)t_2^{n-1}+\sum_{n,n'=1,1}^{2,2}d_{n,n'}t_1^{n-1}t_2^{n'-1},
\end{align}
where $m\in\M([0,1]^2), m_{n',n}\in\M([0,1])$, and $d_{n,n'}\in\R$.  Let $\snorm{\cdot}_{\M_{\DDD}}$ be the seminorm defined in \eqref{eq:5.4.77}. If Assumptions 1 to 6 are verified, then
\begin{equation}
     \mathcal{V}=\underset{f\in\M_{\DDD}([0,1]^2)}{\text{argmin}}\Big(E(\bb{y},\langle\ f,\bm{\nu}\rangle)+\lambda\snorm{f}_{\M_{\DDD}}\Big)
\end{equation}
is nonempty, $\mathrm{weak}^{\star}$-compact, and is the $\mathrm{weak}^{\star}$-closed convex of its extreme points which are of the form
\begin{align}
    f(t_1,t_2)=&\sum_{k=1}^Ka_{k}(t_1-x_{1,k})_+(t_2-x_{2,k})_++\sum_{n=1}^{2}\sum_{m=1}^{K_{1,n}}b_{n,m}t_1^{n-1}(t_2-y_{n,m})_+\nonumber\\
    &+\sum_{n'=1}^{2}\sum_{m'=1}^{K_{2,n'}}c_{n',m'}(t_1-z_{n',m'})_+t_2^{n'-1}+\sum_{n,n'=1,1}^{2,2}d_{n,n'}t_1^{n-1} t_2^{n'-1},
    \end{align}
with the number of knots $\left(K+\sum_{n=1}^{2}K_{1,n}+\sum_{n'=1}^{2}K_{2,n'}\right)$ being upper-bounded by $\left( M-4\right)$.

%%%%%%%%%%%%%%%%%%%%%%%%%%%%%%%%%%%%%%%%%%%%%%%%%%%%%
\section{Multidimensional Theory}
\label{sec:7.1}
We generalize our construction in the $D$-dimensional setting for $\bb{K}=[0,1]^D\subset\R^D$ and the regularization operator
\begin{equation}
    \Dd^{\bb{N}}=\underbrace{\Dd^N\otimes\cdots
    \otimes\Dd^N}_{D\text{ times}},
\end{equation}
with the multiindex notation $\bb{N}=(N,...,N)$. We first define the functional 
\begin{equation}
    \snorm{\cdot}_{\M_{\Dd^{\bb{N}}}}:=\sum_{\underset{\text{max}(\bb{n})=N}{\bb{n}\in[0\ldots N]^D}}\norm{\Dd^{\bb{n}}\{\cdot\}}_{\M\left(\Pi_{d=1}^D[0,\mathbbm{1}_{n_d=N}]\right)},
\end{equation}
with $\bb{n}=(n_1,...,n_D)$.
The interpretation is that $\norm{\cdot}_{\M\left(\Pi_{d=1}^D[0,\mathbbm{1}_{n_d=N}]\right)}$ integrates on a sub-manifold, for example a face or an edge, of Hausdorff dimension equal to the number of those indices $n_d$ that are equal to $N$. The $[0,1]^D$-localized search space $\M_{\Dd^{\bb{N}}}([0,1]^D)$ is defined as the Banach space
\begin{equation}
    \left\{f\in\D'(\R^D):\quad\underset{\text{max}(\bb{n})=N}{\forall\bb{n}\in[0\ldots N]^D,}\quad\Dd^{\bb{n}}\{f\}\in\M\left(\Pi_{d=1}^D[0,\mathbbm{1}_{n_d=N}]\right)\right\},
\end{equation}
on which $\snorm{\cdot}_{\M_{\Dd^{\bb{N}}}}$ is a seminorm with null space 
\begin{equation}
    \mathcal{N}_{\Dd^{\bb{N}}}=\bigotimes_{d=1}^{D}\mathcal{N}_{\Dd^{N}}=\bigotimes_{d=1}^D\text{span}\left\{\frac{t_d^n}{n!}\right\}_{n=0}^{N-1}
\end{equation}
of dimension $N^D.$ Then, the solution set
\begin{equation}
    \mathcal{V}=\underset{f\in\M_{\Dd^{\bb{N}}([0,1]^D)}}{\text{argmin}}\Big(E(\bb{y},\langle f,\bm{\nu}\rangle)+\lambda\norm{f}_{\M_{\Dd^{\bb{N}}}}\Big)
\end{equation}
is nonempty, $\w$-compact, and is the $\w$-closed convex hull of its extreme points which are of the form
\begin{equation}
    f(\cdot)=p(\cdot)+\sum_{k=1}^{K}a_{k}\bigotimes_{d=1}^D\zeta_{d,k}(\cdot-x_{d,k}),\quad p\in\mathcal{N}_{\Dd^{\bb{N}}},
\end{equation}
where $\zeta_{d,k}\in\left\{\frac{(\cdot)^n}{n!}\right\}_{n=0}^{N-1}\cup\left\{\frac{(\cdot)_+^{N-1}}{(N-1)!}\right\}$. If $\zeta_{d,k}\in\left\{\frac{(\cdot)^n}{n!}\right\}_{n=0}^{N-1}$, then the associated knot $x_{d,k}$ is constrained to take value $0$. In addition, the number $K$ of knots is upper-bounded by $\left(M-N^D\right).$

%%%%%%%%%%%%%%%%%%%%%%%%%%%%%%%%%%%%%%%%%%%%%%%%%%%%%
\section{Conclusion}
We showed that it is possible to promote the tensor-product of two one-dimensional splines as the extreme-point solutions of a regularized optimisation problem. In order to have a well-posed problem, there are different ways to regularize on the null space of $\LL$. We uncovered that, for localized problems, there is a unique (causal) null space regularization that takes the form of partial derivatives integrated on the boundary of the image. Our construction is a theoretical framework that: i) justifies the use of tensor-product splines for the resolution of multidimensional inverse problems; ii) derives the unique good regularization;  iii) builds a new bridge between the theory of splines and the regularization theory of multidimensional inverse problems. 

Further work needs to be carried through this bridge to understand how one-dimensional spline-related results transfer to this new multidimensional regime.

%%%%%%%%%%%%%%%%%%%%%%%%%%%%%%%%%%%%%%%%%%%%%%%%%%%%%
\section{Acknowledgement}
Vincent Guillemet was supported by the Swiss National Science Foundation (SNSF) under Grant 200020\_219356.

%%%%%%%%%%%%%%%%%%%%%%%%%%%%%%%%%%%%%%%%%%%%%%%%%%%%%
\appendix
\section{Proofs}

\subsection{Proofs of Section 2.3.}
\label{app:2.3}

\begin{proof}[\textbf{Proof of Theorem \ref{th:3.1}}]
The fact that $\norm{\cdot}_{\epsilon}=\norm{\cdot}_{\infty(\R^2)}$ follows from Theorem \ref{th:densityC} Item 2. Since $\M(\R)\otimes\M(\R)$ is a subspace of $\M(\R^2)$, one finds that  $\M(\R)\widehat{\otimes}_w\M(\R)\subset\M(\R^2)$. Since compactly supported measures are $\w$-dense in $\M(\R^2)$,  we are left to show the inclusion
\begin{equation}
    \M(\bb{K})\subset\M(\R)\widehat{\otimes}_w\M(\R).
\end{equation}
In order to show the equivalent condition that $\forall m\in\M(\bb{K})$,
\begin{equation}
    \exists\{m_n\}_{n=1}^{\infty}\in\M(\R)\otimes\M(\R):\quad\Big(\forall w\in\C_0(\R^2):\quad\underset{n\to\infty}{\text{lim}}\langle m_n,w\rangle\to\langle m,w\rangle\Big),
\end{equation}
we smooth the measure $m$ and then split it into a tensor product. Consider a net $\{\phi_{\epsilon}\}_{\epsilon}$ of Friedrich mollifiers such that the mollifiers $\phi_{\epsilon}$
\begin{itemize}
    \item satisfies $\phi_{\epsilon}\in\C_c^{\infty}(\R^2)\quad$ with $\mathrm{supp}(\phi_{\epsilon})\subset B_{\epsilon}$, where $B_{\epsilon}$ is a ball of radius $\epsilon$;
    \item is normalized absolutely as $\int_{\R^2}\vert\phi_{\epsilon}(\bb{x})\vert\mathrm{d}\bb{x}=1;$
    \item is even-symmetric, with $\forall \bb{y}\in\R^2:\quad\phi_{\epsilon}(\bb{y})=\phi_{\epsilon}(-\bb{y}).$
\end{itemize}
Observe that $\forall w\in\C_0(\R^2)$
\begin{align}
\norm{w-w*\phi_{\epsilon}}_{\infty}&=\underset{\bb{y}\in\R^2}{\text{sup }}\left\vert w(\bb{y})-\int_{\R^2}\phi_{\epsilon}(\bb{y}-\bb{x})w(\bb{x})\mathrm{d}\bb{x}\right\vert\nonumber\\
&\leq\underset{\bb{y}\in\R^2}{\text{sup }}\int_{\R^2}\left\vert\phi_{\epsilon}(\bb{y}-\bb{x})(w(\bb{x})-w(\bb{y}))\right\vert\mathrm{d}\bb{x}\nonumber\\
&=\underset{\bb{y}\in\R^2}{\text{sup }}\int_{\bb{y}+B_{\epsilon}}\left\vert\phi_{\epsilon}(\bb{y}-\bb{x})(w(\bb{x})-f(\bb{y}))\right\vert\mathrm{d}\bb{x}\nonumber\\
&\leq\underset{\bb{x},\bb{y}\in\R^2:\norm{\bb{x}-\bb{y}}_2\leq\epsilon}{\text{sup}}\vert w(\bb{x})-w(\bb{y})\vert\underset{\epsilon\to 0}{\to}0
\end{align}
because $w\in\C_0(\R^2)$ is uniformly continuous. We take $\epsilon_{\ell}$ small enough such that  $\norm{w-w\ast\phi_{\epsilon_{\ell}}}_{\infty}\leq\frac{1}{\ell}$ and  rename $\phi_{\epsilon_{\ell}}$ as $\phi_{\ell}.$ We define the smoothed measure $m_{\ell}$ as
\begin{equation}
    m_{\ell}(\cdot)=(m\ast\phi_{\ell})(\cdot)=\int_{\R^2}\phi_{\ell}(\cdot-\bb{x})\mathrm{d}m(\bb{x}).
\end{equation}
Using the Lebesgue dominated convergence, one finds that $m_{\ell}\in\C^{\infty}(\R^2)$ and, since both $m$ and $\phi_{\ell}$ are compactly supported, $m\ast\phi_{\ell}\in\C(R^2)$ for some square $R^2$. From Theorem \ref{th:densityC} Item 1, we find $w_{\ell}\in\C(R)\otimes\C(R)$ such that $\norm{m_{\ell}-w_{\ell}}_{\infty(R^2)}\leq\frac{1}{\ell}$
and we extend $w_{\ell}$ by 0 outside $R^2$. Consequently, $\forall w\in C_0(\R^2),$ we have that
\begin{align}
\vert\langle m-w_{\ell},w\rangle\vert&\leq\vert\langle m-m_{\ell},w\rangle\vert+\vert\langle m_{\ell}-w_{\ell},w\rangle\vert\nonumber\\
&=\vert\langle m,w-w\ast\phi_{\ell}\rangle\vert+\frac{\vert R\vert^2}{\ell}\norm{w}_{\infty(\R^2)}\nonumber\\
&\leq\frac{1}{\ell}\left(\norm{m}_{\M(\R^2)}+\vert R\vert^2\norm{w}_{\infty(\R^2)}\right).
\end{align}
Since $w_{\ell}\in\M(\R)\otimes\M(\R)$, this proves the claim. 

\end{proof}

\subsection{Proofs of Section 5.1.}
\label{app:4.1} 

\begin{proof}[\textbf{Proof of Proposition \ref{prop:isom}.}]
We only provide a proof of the construction of $(i)$, as the others are similar and simpler. We know from Theorem \ref{Th:predual1D} that, for $i\in\{1,2\}$, the mapping $\Ls_i^{\star}$ is an isometry whose inverse is given by $\Ls_{\bm{\phi}_i}^{-\star}$, with
\begin{equation}
    \Ls_i^{\star}:\left(\C_0(\R),\norm{\cdot}_{\infty(\R)}\right)\to\left(\Ls_i^{\star}\{\C_0(\R)\},\norm{\Ls_{\bm{\phi}_i}^{-\star}\{\cdot\}}_{\infty(\R)}\right).
\end{equation}
It follows from Corollary \ref{cor:isomtensorproduct} (which we leave proofless as a direct extension of \citep[Proposition 3.2]{ryan2002introduction}) that the mapping $\Ls_1^{\star}\otimes\Ls_2^{\star}$ is also an isometry whose inverse is given by $\Ls_{\bm{\phi}_1}^{-\star}\otimes \Ls_{\bm{\phi}_2}^{-\star}$, with
\begin{equation}  \label{eq:4.1.36}
\Ls_1^{\star}\otimes\Ls_2^{\star}:\left(\C_0(\R)\widehat{\otimes}_{\epsilon}\C_0(\R),\norm{\cdot}_{\epsilon(\C_0,\C_0)}\right)\to\left(\Ls_1^{\star}\{\C_0(\R)\}\widehat{\otimes}_{\epsilon}\Ls_2^{\star}\{\C_0(\R)\},\norm{\cdot}_{\epsilon(\Ls_1^{\star}\C_0,\Ls_2^{\star}\C_0)}\right).
\end{equation}

    \begin{Corollary}
    \label{cor:isomtensorproduct}
    Let $\mathrm{O}_1:\mathcal{B}_1\to\mathcal{Z}_1$ and $\mathrm{O}_2:\mathcal{B}_2\to\mathcal{Z}_2$ be two isometries between normed vector spaces. Then, there exists a unique isometry $\mathrm{O}_1\otimes\mathrm{O}_2:\mathcal{B}_1\widehat{\otimes}_{\epsilon}\mathcal{B}_2\to\mathcal{Z}_1\widehat{\otimes}_{\epsilon}\mathcal{Z}_2$ such that 
    \begin{itemize}
        \item [1.] for all $(v_1,v_2)$ in $\mathcal{B}_1\times\mathcal{B}_2:\quad$ $[\mathrm{O}_1\otimes\mathrm{O}_2]\{v_1\otimes v_2\}=\mathrm{O}_1\{v_1\}\otimes\mathrm{O}_2\{v_2\};$
        \item [2.] for all $(z_1,z_2)$ in $\mathcal{Z}_1\times\mathcal{Z}_2:\quad$ $[\mathrm{O}_1^{-1}\otimes\mathrm{O}_2^{-1}]\{z_1\otimes z_2\}=\mathrm{O}_1^{-1}\{z_1\}\otimes\mathrm{O}_2^{-1}\{z_2\};$
        \item [3.] for all $w$ in $\mathcal{B}_1\widehat{\otimes}_{\epsilon}\mathcal{B}_2:\quad\norm{[\mathrm{O}_1\otimes\mathrm{O}_2]\{w\}}_{\epsilon(\mathcal{Z}_1,\mathcal{Z}_2)}=\norm{w}_{\epsilon(\mathcal{B}_1,\mathcal{B}_2)}$;
        \item [4.] for all $w$ in $\mathcal{B}_1\widehat{\otimes}_{\epsilon}\mathcal{B}_2:\quad[\mathrm{O}_1^{-1}\otimes\mathrm{O}_2^{-1}]\{[\mathrm{O}_1\otimes\mathrm{O}_2]\{w\}\}=w.$
    \end{itemize}
\end{Corollary}
In turn, using Theorem \ref{th:densityC}, Equation \eqref{eq:4.1.36} can be rewritten as
\begin{equation}
    \label{eq:4.1.37}\Ls_1^{\star}\otimes\Ls_2^{\star}:\left(\C_0(\R^2),\norm{\cdot}_{\infty(\R^2)}\right)\to\left(\Ls_1^{\star}\{\C_0(\R)\}\widehat{\otimes}_{\epsilon}\Ls_2^{\star}\{\C_0(\R)\},\norm{\cdot}_{\epsilon(\Ls_1^{\star}\C_0,\Ls_2^{\star}\C_0)}\right).
\end{equation}
Furthermore, the same argumentation as in the start of the proof of Item 2 of Theorem \ref{th:3.1} shows that, $\forall v\in\Ls_1^{\star}\{\C_0(\R)\}\otimes\Ls_2^{\star}\{\C_0(\R)\}$,one has that 
\begin{align}
    \norm{\left[\Ls_{\bm{\phi}_1}^{-\star}\otimes \Ls_{\bm{\phi}_2}^{-\star}\right]\{v\}}_{\infty(\R^2)}&=    \norm{\left[\Ls_{\bm{\phi}_1}^{-\star}\otimes \Ls_{\bm{\phi}_2}^{-\star}\right]\{v\}}_{\epsilon(\C_0,\C_0)}\nonumber\\
    &=\norm{[\Ls_1^{\star}\otimes\Ls_2^{\star}]\left\{\left[\Ls_{\bm{\phi}_1}^{-\star}\otimes \Ls_{\bm{\phi}_2}^{-\star}\right]\{[\Ls_1^{\star}\otimes\Ls_2^{\star}]\{f\}\}\right\}}_{\epsilon(\Ls_1^{\star}\C_0,\Ls_2^{\star}\C_0)}\nonumber\\
    &=\norm{[\Ls_1^{\star}\otimes\Ls_2^{\star}]\{f\}}_{\epsilon(\Ls_1^{\star}\C_0,\Ls_2^{\star}\C_0)}=\norm{v}_{\epsilon(\Ls_1^{\star}\C_0,\Ls_2^{\star}\C_0)},
\end{align}
with $v=[\Ls_1^{\star}\otimes\Ls_2^{\star}]\{f\}$. Consequently, the injective completion of $\Ls_1^{\star}\{\C_0(\R)\}\otimes\Ls_2^{\star}\{\C_0(\R)\}$ is equal to the completion under the norm   $\norm{\left[\Ls_{\bm{\phi}_1}^{-\star}\otimes \Ls_{\bm{\phi}_2}^{-\star}\right]\{\cdot\}}_{\infty(\R^2)}$. We denote this completion by $\Ls_1^{\star}\{\C_0(\R)\}\widehat{\otimes}\Ls_2^{\star}\{\C_0(\R)\}$ and \eqref{eq:4.1.37} becomes in
\begin{equation}
\Ls_1^{\star}\otimes\Ls_2^{\star}:\left(\C_0(\R^2),\norm{\cdot}_{\infty(\R^2)}\right)\to\left(\Ls_1^{\star}\{\C_0(\R)\}\widehat{\otimes}\Ls_2^{\star}\{\C_0(\R)\},\norm{\left[\Ls_{\bm{\phi}_1}^{-\star}\otimes \Ls_{\bm{\phi}_2}^{-\star}\right]\{\cdot\}}_{\infty(\R^2)}\right).
\end{equation}
Finally, the space $\Ls_1^{\star}\{\C_0(\R)\}\widehat{\otimes}\Ls_2^{\star}\{\C_0(\R)\}$ is exactly the range of the operator $\Ls_1^{\star}\otimes\Ls_2^{\star}$ over the domain $\C_0(\R^2)$. Therefore, it is equal to $\left[\Ls_1^{\star}\otimes\Ls_2^{\star}\right]\{\C_0(\R^2)\}$. 
\end{proof}

\begin{proof}[\textbf{Proof of Theorem \ref{th:predual}.}] Let the dimension index be $i\in\{1,2\}$ and denote by $i^c=(3-i)$ the complementary index. Consider an arbitrary $g\in\C_{\Ls_1\otimes\Ls_2}(\R^2)$ as in \eqref{eq:4.1.43}.  First, we calculate by linearity that 
\begin{align}
    \left[\Ls_{\bm{\phi}_1}^{-\star}\otimes \Ls_{\bm{\phi}_2}^{-\star}\right]\{g\}=&v+\sum_{n=1}^{N_1}\Ls_{\bm{\phi}_1}^{-\star}\left\{\phi_{1,n}\right\}\otimes \Ls_{\bm{\phi}_2}^{-\star}\left\{\Ls_2^{\star}\{v_{1,n}\}\right\}\nonumber\\
    &+\sum_{n=1}^{N_2}\Ls_{\bm{\phi}_1}^{-\star}\left\{\Ls_1^{\star}\{v_{2,n}\}\right\} \otimes \Ls_{\bm{\phi}_2}^{-\star}\left\{\phi_{2,n}\right\}\nonumber\\
     &+\sum_{n,n'=1,1}^{N_1,N_2}d_{n,n'}\Ls_{\bm{\phi}_1}^{-\star}\{\phi_{1,n}\} \otimes \Ls_{\bm{\phi}_2}^{-\star}\{\phi_{2,n'}\}\nonumber\\
    =&v+\sum_{n=1}^{N_1}0\otimes v_{2,n}
    +\sum_{n=1}^{N_2}v_{1,n}\otimes 0+\sum_{n,n'=1,1}^{N_1,N_2}d_{n,n'}0\otimes0\nonumber\\
    =&v,\label{eq:4.1.52}
\end{align}
where we used \eqref{eq:3.43bis},\eqref{eq:3.43}, and the fact that
 $[\Ls_{\bm{\phi}_1}^{-\star}\otimes \Ls_{\bm{\phi}_2}^{-\star}]\left\{[\Ls_1^{\star}\otimes\Ls_2^{\star}]\{v\}\right\}=v$.
Second, we calculate by linearity that
\begin{align} 
    \left[\mathrm{proj}_{\mathcal{N}_{\Ls_1}^{\star}}\otimes \Ls_{\bm{\phi}_2}^{-\star}\right]\{g\}=&0+\sum_{n=1}^{N_1}\mathrm{proj}_{\mathcal{N}_{\Ls_1}^{\star}}\{\phi_{1,n}\} \otimes \Ls_{\bm{\phi}_2}^{-\star}\{\Ls_2^{\star}\{v_{1,n}\}\}\nonumber\\
    &+\sum_{n=1}^{N_2}\mathrm{proj}_{\mathcal{N}_{\Ls_1}^{\star}}\{\Ls_1^{\star}\{v_{2,n}\}\} \otimes \Ls_{\bm{\phi}_2}^{-\star}\{\phi_{2,n}\}\nonumber\\
    & +\sum_{n,n'=1,1}^{N_1,N_2}d_{n,n'}\mathrm{proj}_{\mathcal{N}_{\Ls_1}^{\star}}\{\phi_{1,n}\}\otimes \Ls_{\bm{\phi}_2}^{-\star}\{\phi_{2,n'}\}\nonumber\\
    =&0+\sum_{n=1}^{N_1}\phi_{1,n}\otimes v_{1,n}+\sum_{n=1}^{N_2}0\otimes 0+\sum_{n,n'=1,1}^{N_1,N_2}d_{n,n'}\phi_{1,n}\otimes0\nonumber\\
    =&\sum_{n=1}^{N_1}\phi_{1,n}\otimes v_{1,n},\label{eq:4.1.57}
\end{align}
where we used \eqref{eq:3.43}, \eqref{eq:3.43bis},\eqref{eq:3.44}, and the (not yet proven) fact that
\begin{equation}
\label{eq:A.2.23}
\forall v\in\C_0(\R^2),\quad\left[\mathrm{proj}_{\mathcal{N}_{\Ls_1}^{\star}}\otimes \Ls_{\bm{\phi}_2}^{-\star}\right]\{[\Ls_1^{\star}\otimes\Ls_2^{\star}]\{v\}\}=0.
\end{equation}
To prove \eqref{eq:A.2.23}, observe that, for  $v=\sum_{n=1}^N\Ls_1^{\star}\{v_{1,n}\}\otimes \Ls_2^{\star}\{v_{2,n}\}$, one has that
\begin{align}
   \left[\mathrm{proj}_{\mathcal{N}_{\Ls_1}^{\star}}\otimes \Ls_{\bm{\phi}_2}^{-\star}\right]\{v\}=\sum_{n=1}^N \mathrm{proj}_{\mathcal{N}_{\Ls_1}^{\star}}\{\Ls_1^{\star}\{v_{1,n}\}\}\otimes \Ls_{\bm{\phi}_2}^{-\star}\{\Ls_2 ^{\star}\{v_{2,n}\}\}=\sum_{n=1}^N0\otimes v_{2,n}=0.
\end{align}
Therefore, $\mathrm{proj}_{\mathcal{N}_{\Ls_1}^{\star}}\otimes \Ls_{\bm{\phi}_2}^{-\star}$ is the 0 mapping over $\Ls_1^{\star}\{\C_0(\R)\}\otimes\Ls_2^{\star}\{\C_0(\R)\}$ and we conclude, with Theorem \ref{th:extension}, that \eqref{eq:A.2.23} holds.
\begin{Theorem}{\citep[Theorem I.7]{reed1972methods}}
\label{th:extension}
Let $G$ be a bounded linear transformation from the normed space $(\mathcal{X}, \norm{\cdot}_{\mathcal{X}})$ to the complete normed space $(\mathcal{Y}, \norm{\cdot}_{\mathcal{Y}})$. Then, G
has a unique extension to a bounded linear map from the 
completion of $\mathcal{X}$ to $\mathcal{Y}$.
\end{Theorem}
 Calculations similar to \eqref{eq:4.1.52} and \eqref{eq:4.1.57} finally yield that
\begin{align}
   \left[\Ls_{\bm{\phi}_1}^{-\star} \otimes\mathrm{proj}_{\mathcal{N}_{\Ls_2}^{\star}}\right]\{v\}&=\sum_{n=1}^{N_2}v_{2,n}\otimes\phi_{2,n},\nonumber\\
   \left[\mathrm{proj}_{\mathcal{N}_{\Ls_1}^{\star}}\otimes\mathrm{proj}_{\mathcal{N}_{\Ls_2}^{\star}}\right]\{v\}&=\sum_{n,n'=1}^{N_1,N_2}d_{n,n'}\phi_{1,n}\otimes\phi_{2,n'}\label{eq:4.1.59}.
\end{align}
The concatenation of \eqref{eq:4.1.52}, \eqref{eq:4.1.57}, and \eqref{eq:4.1.59} yields that
$\C_{\Ls_1\otimes\Ls_2}(\R^2)$ and 
\begin{equation}
     \left([\Ls_1^{\star}\otimes\Ls_2^{\star}]\{\C_0(\R^2)\}\right)\times\left(\mathcal{N}_{\Ls_1}^{\star}\otimes\Ls_2^{\star}\{\C_0(\R)\}\right)\times\left(\Ls_1^{\star}\{\C_0(\R)\}\otimes\mathcal{N}_{\Ls_2}^{\star}\right)\times\left(\mathcal{N}_{\Ls_1}^{\star}\otimes\mathcal{N}_{\Ls_2}^{\star}\right)
\end{equation}
are isometric. Consequently, the sum in \eqref{eq:4.1.35} is direct and the Banach-space structure, together with the norm calculation, follow directly from Proposition \ref{prop:isom}.
\end{proof}

\begin{proof}[\textbf{Proof of Corollary \ref{coro:admissible}}.]
Let $g\in\mathcal{E}'(\R^2)$. Our Remark 1 is that Proposition \ref{prop:compactsupport1D} and Representation \eqref{eq:2.8} yield that $\left[\Ls_{\bm{\phi}_1}^{-\star}\otimes\Ls_{\bm{\phi}_2}^{-\star}\right]\{g\}\in\mathcal{E}'(\R^2)$, $\left[\mathrm{proj}_{\mathcal{N}_{\Ls_1}^{\star}}\otimes\Ls_{\bm{\phi}_2}^{-\star}\right]\{g\}_n\in\mathcal{E}'(\R),$\\
$\left[\Ls_{\bm{\phi}_1}^{-\star}\otimes\mathrm{proj}_{\mathcal{N}_{\Ls_2}^{\star}}\right]\{g\}_n\in\mathcal{E}'(\R),$ and $\left[\mathrm{proj}_{\mathcal{N}_{\Ls_1}^{\star}}\otimes\mathrm{proj}_{\mathcal{N}_{\Ls_2}^{\star}}\right]\{g\}_{n,n'}\in\R$. Then, Theorem \ref{th:predual}  yields that $g\in\C_{\LL}(\R^2)$ if and only if
\begin{itemize}
    \item [1'.] the distribution $\left[\Ls_{\bm{\phi}_1}^{-\star}\otimes\Ls_{\bm{\phi}_2}^{-\star}\right]\{g\}$ is in $\mathcal{C}_0(\R^2)$;
    \item [2'.] the vector of distributions $\left[\mathrm{proj}_{\mathcal{N}_{\Ls_1}^{\star}}\otimes\Ls_{\bm{\phi}_2}^{-\star}\right]\{g\}$ is in $\C_0(\R)^{N_1}$;
    \item [3'.] the vector of distributions $\left[\Ls_{\bm{\phi}_1}^{-\star}\otimes\mathrm{proj}_{\mathcal{N}_{\Ls_2}^{\star}}\right]\{g\}$ is in $\C_0(\R)^{N_2}$;
    \item [4'.] the matrix $\left[\mathrm{proj}_{\mathcal{N}_{\Ls_1}^{\star}}\otimes\mathrm{proj}_{\mathcal{N}_{\Ls_2}^{\star}}\right]\{g\}$ is in $\R^{N_1\times N_2}$.
\end{itemize}
It follows from Remark 1 that Conditions $1',2',3',$ and $4'$ simplify into $1,2$, and $3$ of the proposition statement. Finally, Items 4 and 5 are reformulations of Items 1,2, and 3 in the case where $g$ is separable.
\end{proof}

\begin{proof}[\textbf{Proof of Corollary \ref{coro:testfunctions}}]
\hfill\\
\textbf{Item 1.} Let $\psi\in\D(\R^2)$. It is sufficient to prove that Items 1, 2, and 3 of Corollary \ref{coro:admissible} are verified. The verification of Items 2 and 3 being easy, we focus here on the verification of Item 1. We calculate that

\begin{align}
\Big[\Ls_{\bm{\phi}_1}^{-\star}&\otimes \Ls_{\bm{\phi}_2}^{-\star}\Big]\{\psi\}(x_1,x_2)\nonumber\\
=&\underbrace{\int_{\R^2}g_{\Ls_1}(t_1-x_1)g_{\Ls_2}(t_2-x_2)\psi(t_1,t_2)\mathrm{d}t_1\mathrm{d}t_2}_{(A)}\nonumber\\
&-\sum_{n=1}^{N_2}\underbrace{(g_{\Ls_2}^{\vee}\ast\phi_{2,n})(x_2)\int_{\R^2}g_{\Ls_1}(t_1-x_1)p_{2,n}(t_2)\psi(t_1,t_2)\mathrm{d}t_1\mathrm{d}t_2}_{(B)}\nonumber\\
&-\sum_{n=1}^{N_1}\underbrace{(g_{\Ls_1}^{\vee}\ast\phi_{1,n})(x_1)\int_{\R^2}g_{\Ls_2}(t_2-x_2)p_{1,n}(t_1)\psi(t_1,t_2)\mathrm{d}t_1\mathrm{d}t_2}_{(C)}\nonumber\\
&+\sum_{n,n'=1,1}^{N_1,N_2}\underbrace{(g_{\Ls_1}^{\vee}\ast\phi_{1,n})(x_1)(g_{\Ls_2}^{\vee}\ast\phi_{2,n})(x_2)\int_{\R^2}p_{1,n}(t_1)p_{2,n}(t_2)\psi(t_1,t_2)\mathrm{d}t_1\mathrm{d}t_2}_{(D)}\label{eq:A.3.111}.
\end{align}
The term $(A)$ is a convolution between the distribution $g_{\Ls_1}^{\vee}\otimes g_{\Ls_2}^{\vee}$ and the test function $\psi$; thus, $(A)\in\C^{\infty}(\R^2)$. Likewise, $(D)\in\C^{\infty}(\R^2)$. The same argumentation is conducted with $(B),(C)$ to conclude that $((B),(C))\in\C^{\infty}(\R^2)^2$, with the additional remarks that 
\begin{equation}
\int_{\R}p_{1,n}
(t_1)\psi(t_1,t_2)\mathrm{d}t_1\in\D(\R)\quad\text{ and}\quad\int_{\R}p_{2,n}(t_2)\psi(t_1,t_2)\mathrm{d}t_2\in\D(\R).
\end{equation}
Therefore, $\left[\Ls_{\bm{\phi}_1}^{-1}\otimes \Ls_{\bm{\phi}_2}^{-1}\right]\{\psi\}\in\C^{\infty}(\R^2)\subset\C(\R^2).$ 
\hfill\\
\textbf{Item 2.} The argumentation is the same as in the proof of Item 1, up to the fact that the convolution between a distribution and a test function outputs a $\C^{\infty}$ function. In Item 2, we have to use the fact that the convolution between an $\mathcal{L}_{\infty,\mathrm{loc}}$ function and an $\mathcal{L}_{p,c}$ function outputs a continuous function. 
\hfill\\
\textbf{Item 3.} Observe that $\delta_{\bb{x}}=\delta_{x_1}\otimes\delta_{x_2}$, which is a separable and compactly supported distribution. Items 4 and 5 of Corollary \ref{coro:admissible} are verified with the observation that, if $N_1\geq2$ and $N_2\geq2$, then $g_{\Ls_1}$ and $g_{\Ls_2}$ are continuous functions.
\end{proof}

\subsection{Proofs of Section 5.2.}
\label{app:4.2}

\begin{proof}[\textbf{Proof of Proposition \ref{prop:isom2}.}]Like in the proof of Proposition \ref{prop:isom}, we only prove Item 1. As first step we show that 
\begin{equation}
    [\Ls_{\bm{\phi}_1}^{-1}\otimes \Ls_{\bm{\phi}_2}^{-1}]\{m\}=\int_{\R^2}g_{\bm{\phi}_1}(\cdot,x_1)g_{\bm{\phi}_2}(\cdot,x_2)\mathrm{d}m(x_1,x_2)
\end{equation}
is well-defined. To do so we recall that $g_{\bm{\phi}_1}(t_1,\cdot)g_{\bm{\phi}_2}(t_2,\cdot)$ is compactly supported in 
\begin{equation}
    [\text{min}(t_1,\phi_1^-),\text{max}(t_1,\phi_1^+)]\times[\text{min}(t_2,\phi_2^-),\text{max}(t_2,\phi_2^+)].
\end{equation}
Hence, for $\bb{B}=[a,b]\times[c,d]$ and
\begin{equation}
    \bb{A}=[\text{min}(a,\phi_1^-),\text{max}(b,\phi_1^+)]\times[\text{min}(c,\phi_2^-),\text{max}(d,\phi_2^+)],
\end{equation}
one has that 
\begin{equation}
\forall(t_1,t_2)\in\bb{B},\forall(x_1,x_2)\notin\bb{A}:\quad g_{\bm{\phi}_1}(t_1,x_1)g_{\bm{\phi}_2}(t_2,x_2)=0.
\end{equation}
If $\psi\in\C_c^{\infty}(\R^2)$ with $\text{supp}(\psi)\subset\bb{B}$, then 
\begin{align}
    \left\langle[\Ls_{\bm{\phi}_1}^{-1}\otimes \Ls_{\bm{\phi}_2}^{-1}]\{m\},\psi\right\rangle&=\left\langle    \bm{1}_{\bb{B}}[\Ls_{\bm{\phi}_1}^{-1}\otimes \Ls_{\bm{\phi}_2}^{-1}]\{m\},\psi\right\rangle=\left\langle    \bm{1}_{\bb{B}}[\Ls_{\bm{\phi}_1}^{-1}\otimes \Ls_{\bm{\phi}_2}^{-1}]\{m\bm{1}_{\bb{A}}\},\psi\right\rangle\nonumber\\
    &=\left\langle[\Ls_{\bm{\phi}_1}^{-1}\otimes \Ls_{\bm{\phi}_2}^{-1}]\{m\bm{1}_{\bb{A}}\},\psi\right\rangle,
\end{align}
where $[\Ls_{\bm{\phi}_1}^{-1}\otimes \Ls_{\bm{\phi}_2}^{-1}]\{m\bm{1}_{\bb{A}}\}$ is well-defined because, from \eqref{eq:2.8}, it can be expended into a sum of convolutions between Green's functions $(g_{\Ls_1},g_{\Ls_2},g_{\Ls_1}\otimes g_{\Ls_2})$ and $m\bm{1}_A$, which are well-defined as convolutions between right-sided distributions.

As second step, we show the isometry from the proposition statement. For all $\psi\in\C_c^{\infty}(\R^2)$ with $\text{supp}(\psi)\subset\bb{B}$, we calculate that 
\begin{align}
    &\left\langle[\LL]\{[\Ls_{\bm{\phi}_1}^{-1}\otimes \Ls_{\bm{\phi}_2}^{-1}]\{m\}\},\psi\right\rangle=\left\langle[\LL]\{[\Ls_{\bm{\phi}_1}^{-1}\otimes \Ls_{\bm{\phi}_2}^{-1}]\{m\bm{1}_{\bb{A}}\}\},\psi\right\rangle\nonumber\\
    &=\left\langle[\LL]\{(g_{\Ls_1}\otimes g_{\Ls_2})\ast m\bm{1}_{\bb{A}}\},\psi\right\rangle=\langle m\bm{1}_{\bb{A}},\psi\rangle=\langle m,\psi\rangle,
\end{align}
where we used \eqref{eq:2.8} as well as the fact that $\bb{B}\subset\bb{A}$. The density of $\C_c^{\infty}(\R^2)$ in $\C_0(\R^2)$ implies that $[\LL]\{[\Ls_{\bm{\phi}_1}^{-1}\otimes \Ls_{\bm{\phi}_2}^{-1}]\{m\}\}=m$. The isometry follows directly. 
\end{proof}

\begin{proof}[\textbf{Proof of Theorem \ref{th:native}.}]
Let the dimension index be $i\in\{1,2\}$ and denote by $i^c=(3-i)$ the complementary index. Consider an arbitrary $f\in\M_{\Ls_1\otimes\Ls_2}(\R^2)$ as in \eqref{eq:4.2.78}. First, we calculate that 
\begin{align}
    [\Ls_1\otimes\Ls_2]\{f\}=&m+\sum_{n=1}^{N_1}\Ls_1\{p_{1,n}\}\otimes \Ls_2\{\Ls_{\bm{\phi}_2}^{-1}\{m_{1,n}\}\}\nonumber\\
    &+\sum_{n=1}^{N_2}\Ls_1\{\Ls_{\bm{\phi}_1}^{-1}\{m_{2,n}\}\}\otimes \Ls_2\{p_{2,n}\}\nonumber\\
    &+\sum_{n,n'=1,1}^{N_1,N_2}d_{n,n'}\Ls_1\{p_{1,n}\}\otimes \Ls_2\{p_{2,n'}\}\nonumber\\
    =&m+\sum_{n=1}^{N_1}0\otimes m_{1,n}+\sum_{n=1}^{N_2}m_{2,n}\otimes 0+\sum_{n,n'=1,1}^{N_1,N_2}d_{n,n'}0\otimes 0\nonumber\\
    =&m,\label{eq:4.2.93}
\end{align}
where we used the fact that $[\Ls_1\otimes\Ls_2]\{[\Ls_{\bm{\phi}_1}^{-1}\otimes \Ls_{\bm{\phi}_2}^{-1}]\{m\}\}=m$. Second, we calculate that
\begin{align}
   [\mathrm{proj}_{\mathcal{N}_{\Ls_1}}\otimes\Ls_2]\{f\}=&[\mathrm{proj}_{\mathcal{N}_{\Ls_1}}\otimes\Ls_2]\{[\Ls_{\bm{\phi}_1}^{-1}\otimes \Ls_{\bm{\phi}_2}^{-1}]\{m\}\}\nonumber\\
   &+\sum_{n=1}^{N_1} \mathrm{proj}_{\mathcal{N}_{\Ls_1}}\{p_{1,n}\}\otimes \Ls_2\{\Ls_{\bm{\phi}_2}^{-1}\{m_{1,n}\}\}\nonumber\\\
   &+\sum_{n=1}^{N_2}\mathrm{proj}_{\mathcal{N}_{\Ls_1}}\{\Ls_{\bm{\phi}_1}^{-1}\{m_{2,n}\}\} \otimes \Ls_2\{p_{2,n}\}\nonumber\\
   &+\sum_{n,n'=1,1}^{N_1,N_2}d_{n,n'} \mathrm{proj}_{\mathcal{N}_{\Ls_1}}\{ p_{1,n}\}\otimes\Ls_2\{ p_{2,n'}\}\nonumber\\
    =&0+\sum_{n=1}^{N_1}p_{1,n}\otimes m_{1,n}+\sum_{n=1}^{N_2}0\otimes 0+\sum_{n,n'=1,1}^{N_1,N_2}d_{n,n'} p_{1,n}\otimes0\nonumber\\
    =&\sum_{n=1}^{N_1}p_{1,n}\otimes m_{1,n}\label{eq:4.2.100},
\end{align}
where we used \eqref{eq:3.38} and the (not yet proven) fact that 
\begin{equation}
    \forall m\in\M(\R^2),\quad\left[\mathrm{proj}_{\mathcal{N}_{\Ls_1}}\otimes\Ls_2\right]\{[\Ls_{\bm{\phi}_1}^{-1}\otimes \Ls_{\bm{\phi}_2}^{-1}]\{m\}\}=0.
\end{equation}
Observe that, for $m=\sum_{n=1}^Nm_{1,n}\otimes m_{2,n}\in\M(\R)\otimes\M(\R)$, one has that 
\begin{align}
&\left[\mathrm{proj}_{\mathcal{N}_{\Ls_1}}\otimes\Ls_2\right]\left\{m\right\}=\sum_{n=1}^N\mathrm{proj}_{\mathcal{N}_{\Ls_1}}\{\Ls_{\bm{\phi}_1}^{-1}\{m_{1,n}\}\}\otimes \Ls_2\left\{\Ls_{\bm{\phi}_2}^{-1}\{m_{2,n}\}\right\}=\sum_{n=1}^{N}0\otimes m_{2,n}.
\end{align}
Then, we conclude from the $\w$ density of $\M(\R)\otimes\M(\R)$ in $\M(\R^2)$ that the unique continuous extension of  
\begin{equation}
    \left[\mathrm{proj}_{\mathcal{N}_{\Ls_1}}\otimes\Ls_2\right]\left\{\left[\Ls_{\bm{\phi}_1}^{-1}\otimes \Ls_{\bm{\phi}_2}^{-1}\right]\{\cdot\}\right\}:\M(\R)\otimes\M(\R)\to\R^{N_1}\otimes\M(\R)
\end{equation}
to $\M(\R^2)$ is the zero mapping.
Third, we find with the same technique that 
\begin{align}
     \label{eq:4.2.102}
     &\left[\Ls_1\otimes\mathrm{proj}_{\mathcal{N}_{\Ls_2}}\right]\{f\}=\sum_{n=1}^{N_2}m_{2,n}\otimes p_{2,n },\nonumber\\
     &\left[\mathrm{proj}_{\mathcal{N}_{\Ls_1}}\otimes\mathrm{proj}_{\mathcal{N}_{\Ls_2}}\right]\{f\}=\sum_{n,n'=1,1}^{N_1,N_2}d_{n,n'}p_{1,n}\otimes p_{2,n'}.
\end{align}
It follows from \eqref{eq:4.2.93}, \eqref{eq:4.2.100}, and \eqref{eq:4.2.102} that $\M_{\LL}(\R^2)$ is isomorphic to 
\begin{align}
    \left([\Ls_{\bm{\phi}_1}^{-1}\otimes \Ls_{\bm{\phi}_2}^{-1}]\{\M(\R^2)\}\right)\times\left(\mathcal{N}_{\Ls_1}\otimes \Ls_{\bm{\phi}_2}^{-1}\{\M(\R)\}\right)&\times\left(\Ls_{\bm{\phi}_1}^{-1}\{\M(\R)\}\otimes\mathcal{N}_{\Ls_2}\right)\nonumber\\
    &\times\left(\mathcal{N}_{\Ls_1}\otimes\mathcal{N}_{\Ls_2}\right).
\end{align}
This establishes the direct-sum decomposition. This implies, with Proposition \ref{prop:isom2}, the claimed Banach-space structure and the norm calculation.
\end{proof}

\begin{proof}[\textbf{Proof of Proposition \ref{prop:4.4.1.1}.}]
Denote by $\tilde{\M}_{\LL}(\R^2)$ the set defined in the proposition statement. It is clear from Theorem \ref{th:native} and its proof that, if $f\in\M_{\LL}(\R^2)$, then $f\in\tilde{\M}_{\LL}(\R^2)$. To prove the converse inclusion, remark that, for $f\in\D'(\R^2)$, the condition $[\LL]\{f\}\in\M(\R^2)$ implies that 
\begin{equation}
\label{eq:A.5.1.142}
    f=\left[\Ls_{\bm{\phi_1}}^{-1}\otimes \Ls_{\bm{\phi_2}}^{-1}\right]\{m\}+\sum_{n=1}^{N_1}p_{1,n}\otimes f_{1,n}+\sum_{n=1}^{N_2}f_{2,n}\otimes p_{2,n},
\end{equation}
with $m\in\M(\R^2)$ and $(f_{1,n})_{n=1}^{N_1}\in\D'(\R)^{N_1}$, and $(f_{2,n})_{n=1}^{N_2}\in\D'(\R)^{N_2}$. In turn, the condition $[\Ls_1\otimes\mathrm{proj}_{\mathcal{N}_{\Ls_2}}]\{f\}\in\M(\R)\otimes\R^{N_2}$ implies that 
\begin{equation}
\label{eq:A.5.1.143}
    \sum_{n=1}^{N_2}f_{2,n}\otimes p_{2,n}=\sum_{n=1}^{N_2}\Ls_{\bm{\phi_1}}^{-1}\{m_{2,n}\}\otimes p_{2,n}+\sum_{n,n'=1,1}^{N_1,N_2}d_{n,n'}p_{1,n}\otimes p_{2,n'},
\end{equation}
with $(m_{2,n})_{n=1}^{N_2}\in\M(\R)^{N_2}$ and $(d_{n,n'})_{n,n'=1,1}^{N_1,N_2}\in\R^{N_1\times N_2}$. Likewise, the condition $[\mathrm{proj}_{\mathcal{N}_{\Ls_1}}\otimes\Ls_2]\{f\}\in\R^{N_2}\otimes\M(\R)$ implies that 
\begin{equation}
\label{eq:A.5.1.144}
    \sum_{n=1}^{N_1}p_{1,n}\otimes f_{1,n}=\sum_{n=1}^{N_1}p_{1,n}\otimes \Ls_{\bm{\phi_2}}^{-1}\{m_{1,n}\}+\sum_{n,n'=1,1}^{N_1,N_2}\tilde{d}_{n,n'}p_{1,n}\otimes p_{2,n'},
\end{equation}
with $(m_{1,n})_{n=1}^{N_1}\in\M(\R)^{N_1}$ and $(\tilde{d}_{n,n'})_{n,n'=1,1}^{N_1,N_2}\in\R^{N_1\times N_2}$. The combination of \eqref{eq:A.5.1.142}, \eqref{eq:A.5.1.143}, and \eqref{eq:A.5.1.144} implies that $f$ is exactly of the form \eqref{eq:4.2.78}, which concludes the proof.
\end{proof}

\begin{proof}[\textbf{Proof of Proposition \ref{prop:extremepoints}}]
  It is straightforward to calculate that 
\begin{itemize}
    \item [1'.] the elements of the form of Item 1 are exactly the the extreme points of the unit ball in  
    \begin{equation}
        \left(\left[\Ls_{\bm{\phi}_1}^{-1}\otimes \Ls_{\bm{\phi}_2}^{-1}\right]\{\M(\R^2)\},\quad\norm{[\Ls_1\otimes\Ls_2]\{f\}}_{\M(\R^2)}\quad\right);
    \end{equation}
    \item [2'.] the elements of the form of Item 2 are exactly the the extreme points of the unit ball in 
    \begin{equation}
        \left(\mathcal{N}_{\Ls_1}\otimes \Ls_{\bm{\phi}_2}^{-1}\{\M(\R)\},\quad\norm{\left[\mathrm{proj}_{\mathcal{N}_{\Ls_1}}\otimes \Ls_2\right]\{f\}}_{\bb{p}_1,\M(\R)}\right);
    \end{equation}
    \item [3'.]the elements of the form of Item 3 are exactly the the extreme points of the unit ball in  
    \begin{equation}
        \left(\Ls_{\bm{\phi}_1}^{-1}\{\M(\R)\}\otimes\mathcal{N}_{\Ls_2},\quad\norm{ \left[\Ls_1\otimes\mathrm{proj}_{\mathcal{N}_{\Ls_2}}\right]\{f\}}_{\M(\R),\bb{p}_2}\right).
    \end{equation}
\end{itemize}
Finally, since $\snorm{\cdot}_{\M_{\LL}}$ is zero-valued on $\mathcal{N}_{\Ls_1}\otimes\mathcal{N}_{\Ls_2}$, it follows from \citep[Lemma 1]{unser2022convex} and $1',2',$ and $3'$ that the extreme points of the unit ball of $\snorm{\cdot}_{\M_{\LL}}$ have the claimed structure.
\end{proof}

\begin{proof}[\textbf{Proof of Theorem \ref{Th:7}.}]
From the isometry between $[\Ls_1^{\star}\otimes\Ls_2^{\star}]\{\C_0(\R)^2\}$ and $\C_0(\R^2)$, and the duality between $\M(\R^2)$ and $\C_0(\R^2),$ we deduce that a functional $F$ is in the dual space of $(i)$ if and only if it is such that
\begin{equation}
    \forall[\Ls_1^{\star}\otimes\Ls_2^{\star}]\{v\}\in[\Ls_1^{\star}\otimes\Ls_2^{\star}]\{\C_0(\R)\otimes\C_0(\R)\}:\quad\quad\langle F,[\Ls_1^{\star}\otimes\Ls_2^{\star}]\{v\}\rangle=\langle m,v\rangle
\end{equation}
for some $m\in\M(\R^2)$. A short calculation shows that $F=[\Ls_{\bm{\phi}_1}^{-1}\otimes \Ls_{\bm{\phi}_2}^{-1}]\{m\}$ and, therefore, that $(i)^{\star}$ is the dual space of $(i)$.
Next, we observe that
\begin{align}
    &\left(\M(\R),\norm{\cdot}_{\M(\R)}\right)=\left(\C_0(\R),\norm{\cdot}_{\infty(\R)}\right)^{\star},\nonumber\\
    &\left(\mathcal{N}_{\Ls_1},\norm{\cdot}_{\bb{p}_1}\right)=\left(\mathcal{N}_{\Ls_1}^{\star},\norm{\cdot}_{\bm{\phi}_1}\right)^{\star}\nonumber\\
    \Rightarrow&\left(\mathcal{N}_{\Ls_1}\otimes\M(\R),\norm{\cdot}_{\bb{p}_1,\M(\R)}\right)=\left(\mathcal{N}_{\Ls_1}^{\star}\otimes\C_0(\R),\norm{\cdot}_{\bm{\phi}_1,\infty(\R)}\right)^{\star}\label{eq:A.4.112}.
\end{align} 
The claim that $(ii)^{\star}$ is the dual space of $(ii)$ follows from \eqref{eq:A.4.112} and the isometries in Propositions \ref{prop:isom} and \ref{prop:isom2}. The duality of $(iii)^{\star}$ and $(iv)^{\star}$ with $(iii)$ and $(iv)$ are proved likewise. It follows from these dualities that 
\begin{align}
    \Big(&\left[\Ls_{\bm{\phi}_1}^{-1}\otimes \Ls_{\bm{\phi}_2}^{-1}\right]\{\M(\R^2)\}\Big)\times\left(\mathcal{N}_{\Ls_1}\otimes \Ls_{\bm{\phi}_2}^{-1}\{\M(\R)\} \right)\nonumber\\
    &\times\left(\Ls_{\bm{\phi}_1}^{-1}\{\M(\R)\}\otimes\mathcal{N}_{\Ls_2}\right)\times\left(\mathcal{N}_{\Ls_1}\otimes\mathcal{N}_{\Ls_2}
    \right)
\end{align}
equipped with the norm \eqref{eq:4.2.79} is the dual space of 
\begin{equation}
\left([\Ls_1^{\star}\otimes\Ls_2^{\star}]\{\C_0(\R^2)\}\right)\times\left(\mathcal{N}_{\Ls_1}^{\star}\otimes\Ls_2^{\star}\{\C_0(\R)\}\right)\times\left(\Ls_1^{\star}\{\C_0(\R)\}\otimes\mathcal{N}_{\Ls_2}^{\star}\right)\times\left(\mathcal{N}_{\Ls_1}^{\star}\otimes\mathcal{N}_{\Ls_2}^{\star}
    \right)
\end{equation}
equipped with the norm \eqref{eq:4.1.44}. Finally, the desired result follows from using this duality and the direct-sum decompositions \eqref{eq:4.1.43} and \eqref{eq:4.2.78}. 
\end{proof}

\subsection{Proofs of Section 5.3.}
\label{app:4.3}

\begin{proof}[\textbf{Proof of Proposition \ref{prop:regularity}.}] In this proof, we make use of our knowledge of the support of the kernel $g_{\bm{\phi}_1}\otimes g_{\bm{\phi}_2}$ and of its regularity (Proposition \ref{prop:compactsupport1D}). Let $\bb{B}$ and $\bb{A}$ be two sets as in Proposition \ref{prop:isom2}. Then, 
\begin{equation} \forall(t_1,t_2)\in\bb{B},\forall(x_1,x_2)\notin\bb{A}:\quad g_{\bm{\phi}_1}(t_1,x_1)g_{\bm{\phi}_2}(t_2,x_2)=0.
\end{equation}
For all $\psi\in\C_c(\bb{B})$ one has, using \eqref{eq:2.8} and $m_{\bb{A}}=m\bb{1}_{\bb{A}}$, that
\begin{align}
    &\left\langle\frac{\partial^{d_1+d_2}[\Ls_{\bm{\phi}_1}^{-1}\otimes \Ls_{\bm{\phi}_2}^{-1}]\{m\}}{\partial t_1^{d_1}\partial t_2^{d_2}},\psi\right\rangle=\nonumber\\
    &\left\langle\frac{\partial^{d_1+d_2}[\Ls_{\bm{\phi}_1}^{-1}\otimes \Ls_{\bm{\phi}_2}^{-1}]\{m_{\bb{A}}\}}{\partial t_1^{d_1}\partial t_2^{d_2}},\psi\right\rangle=\nonumber\\
    &\langle(\Dd^{d_1}g_{\Ls_1}\otimes\Dd^{d_2}g_{\Ls_2})\ast m_{\bb{A}},\psi\rangle\nonumber\\
    &-\sum_{n=1}^{N_1}\left\langle\Dd^{d_1}p_{1,n}\otimes(\Dd^{d_2}g_{\Ls_2}\ast \langle m_{\bb{A}}(t,\cdot),(g_{\Ls_1}^{\vee}\ast\phi_{1,n})(t)\rangle),\psi\right\rangle\nonumber\\
    &-\sum_{n=1}^{N_1}\left\langle(\Dd^{d_1}g_{\Ls_1}\ast \langle m_{\bb{A}}(\cdot,t),(g_{\Ls_2}^{\vee}\ast\phi_{2,n})(t)\rangle)\otimes\Dd^{d_2}p_{2,n},\psi\right\rangle\nonumber\\
    &+\sum_{n,n'=1,1}^{N_1,N_2}\langle m_{\bb{A}},(g_{\Ls_1}^{\vee}\ast\phi_{1,n})\otimes(g_{\Ls_2}^{\vee}\ast\phi_{2,n})\rangle\left\langle\Dd^{d_1}p_{1,n}\otimes\Dd^{d_2}p_{2,n},\psi\right\rangle.
\end{align}
In addition,
\begin{equation}
    \Dd^{d_i}g_{\Ls_i}\in\mathcal{L}_{\infty,loc}(\R),\quad\Dd^{d_i}p_{i,n}\in\mathcal{L}_{\infty,loc}(\R),
\end{equation}
and 
\begin{equation}
       \langle m_{\bb{A}}(t,\cdot),(g_{\Ls_1}^{\vee}\ast\phi_{1,n})(t)\rangle\in\M(\R),\quad\langle m_{\bb{A}}(\cdot,t),(g_{\Ls_2}^{\vee}\ast\phi_{2,n})(t)\rangle\in\M(\R).
\end{equation}
Consequently, taking the supremum over $\psi$ with $\norm{\psi}_{\mathcal{L}_1(\bb{B})}\leq1$, we obtain that, for any open set $\bb{O}\subset\bb{B}$,
\begin{equation}
    \norm{\frac{\partial^{d_1+d_2}[\Ls_{\bm{\phi}_1}^{-1}\otimes \Ls_{\bm{\phi}_2}^{-1}]\{m\}}{\partial t_1^{d_1}\partial t_2^{d_2}}}_{\mathcal{L}_{\infty}(\bb{O})}<\infty\Rightarrow\frac{\partial^{d_1+d_2}[\Ls_{\bm{\phi}_1}^{-1}\otimes \Ls_{\bm{\phi}_2}^{-1}]\{m\}}{\partial t_1^{d_1}\partial t_2^{d_2}}\in\mathcal{L}_{\infty,\mathrm{loc}}(\R^2).
\end{equation}
The proof for the continuity is based on the Lebesgue dominated convergence and is omitted.
\end{proof}

\subsection{Proofs of Section 5.4.}
\label{app:4.5}

\begin{proof}[\textbf{Proof of Proposition \ref{prop:CompactRepresentation}.}]
It follows from \eqref{eq:4.5.55} that $f\in\M_{\LL}(\bb{K})$ has the representation
\begin{align}
    &\underbrace{\int_{\bb{K}}g_{\bm{\phi}_1}(t_1,x_1)g_{\bm{\phi_2}}(t_2,x_2)\mathrm{d}m(x_1,x_2)}_{A)}\nonumber+\underbrace{\sum_{n=1}^{N_1}p_{1,n}(t_1)\left(\int_{\mathrm{K}_2}g_{\bm{\phi_2}}(t_2,x_2)\mathrm{d}m_{1,n}(x_2)\right)}_{B)}\nonumber\\
    &+\underbrace{\sum_{n'=1}^{N_2}\left(\int_{\mathrm{K}_1}g_{\bm{\phi_1}}(t_1,x_1)\mathrm{d}m_{2,n'}(x_1)\right) p_{2,n'}(t_2)}_{C)}+\underbrace{\sum_{n,n'=1,1}^{N_1,N_2}d_{n,n'}p_{1,n}(t_1)p_{2,n'}(t_2)}_{D)}\label{eq:A.7.150}.
\end{align}
Then, for $i\in\{1,2\}$, we inject the kernel representation \eqref{eq:2.8} inside $A)$ to find that
\begin{align}
    A)=&\int_{\bb{K}}g_{\Ls_1}(t_1-x_1)g_{\Ls_2}(t_2-x_2)\mathrm{d}m(x_1,x_2)\nonumber\\
    &-\sum_{n'=1}^{N_2}p_{2,n'}(t_2)\int_{\mathrm{K}_1}g_{\Ls_1}(t_1-x_1)\mathrm{d}m^{\mathrm{A}}_{2,n'}(x_1)\nonumber\\
    &-\sum_{n=1}^{N_1}p_{1,n}(t_1)\int_{\mathrm{K}_2}g_{\Ls_2}(t_2-x_2)\mathrm{d}m^{\mathrm{A}}_{1,n}(x_2)+\sum_{n,n'=1,1}^{N_1,N_2}d^{\mathrm{A}}_{n,n'}p_{1,n}(t_1)p_{2,n'}(t_2),\label{eq:A.7.96}
\end{align}
with 
\begin{align}
    m^{\mathrm{A}}_{2,n'}(\cdot)&=\int_{\mathrm{K}_2}(g_{\Ls_2}^{\vee}\ast\phi_{2,n'})(x_2)\mathrm{d}m(\cdot,x_2)\in\M(\mathrm{K}_1)\nonumber\\
    m^{\mathrm{A}}_{1,n}(\cdot)&=\int_{\mathrm{K}_1}(g_{\Ls_1}^{\vee}\ast\phi_{1,n})(x_1)\mathrm{d}m(x_1,\cdot)\in\M(\mathrm{K}_2)\nonumber\\
     d^{\mathrm{A}}_{n,n'}&=\int_{\bb{K}}(g_{\Ls_1}^{\vee}\ast\phi_{1,n})(x_1)(g_{\Ls_2}^{\vee}\ast\phi_{2,n})(x_2)\mathrm{d}m(x_1,x_2).
\end{align}
Likewise, we inject the kernel representation \eqref{eq:2.8} inside $B)$ and $C)$ to find that 
\begin{align}
    B)=\sum_{n=1}^{N_1}p_{1,n}(t_1)\left(\int_{\mathrm{K}_2}g_{\Ls_2}(t_2-x_2)\mathrm{d}m_{1,n}(x_2)\right)-\sum_{n,n'=1,1}^{N_1,N_2}d^{\mathrm{B}}_{n,n'}p_{1,n}(t_1)p_{2,n'}(t_2),\label{eq:A.7.154}\\
    C)=\sum_{n'=1}^{N_2}p_{2,n}(t_2)\left(\int_{\mathrm{K}_1}g_{\Ls_1}(t_1-x_1)\mathrm{d}m_{2,n'}(x_1)\right)-\sum_{n,n'=1,1}^{N_1,N_2}d^{\mathrm{C}}_{n,n'}p_{1,n}(t_1)p_{2,n'}(t_2),\label{eq:A.7.155}
\end{align}
with 
\begin{align}
    d_{n,n'}^{\mathrm{B}}&=\int_{\mathrm{K}_2}(g_{\Ls_2}^{\vee}\ast\phi_{2,n'})(x_2)\mathrm{d}m_{1,n}(x_2)\in\R\nonumber\\
    d_{n,n'}^{\mathrm{C}}&=\int_{\mathrm{K}_1}(g_{\Ls_1}^{\vee}\ast\phi_{1,n})(x_1)\mathrm{d}m_{2,n'}(x_1)\in\R.
\end{align}
Finally, the combination of \eqref{eq:A.7.150}, \eqref{eq:A.7.96}, \eqref{eq:A.7.154}, and \eqref{eq:A.7.155} yield the awaited integral representation
\begin{eqnarray}
    f(t_1,t_2)&=&\int_{\bb{K}}g_{\Ls_1}(t_1-x_1)g_{\Ls_2}(t_2-x_2)\mathrm{d}m(x_1,x_2)\nonumber\\
    &&\mbox{}+\sum_{n'=1}^{N_2}p_{2,n'}(t_2)\int_{\mathrm{K}_1}g_{\Ls_1}(t_1-x_1)\mathrm{d}(m_{2,n'}-m^{\mathrm{A}}_{2,n'})(x_1)\nonumber\\
    &&\mbox{}+\sum_{n=1}^{N_1}p_{1,n}(t_1)\int_{\mathrm{K}_2}g_{\Ls_2}(t_2-x_2)\mathrm{d}(m_{1,n}-m^{\mathrm{A}}_{1,n})(x_2)\nonumber\\
    &&\mbox{}+\sum_{n,n'=1,1}^{Nm_1,N_2}(d^{\mathrm{A}}_{n,n'}-d^{\mathrm{B}}_{n,n'}-d^{\mathrm{C}}_{n,n'}+d_{n,n'})p_{1,n}(t_1)p_{2,n'}(t_2)\nonumber\\.
\end{eqnarray}
\end{proof}

\begin{proof}[\textbf{Proof of Proposition \ref{prop:normequivalence}.}]
Recall that $\mathrm{proj}_{\mathcal{N}_{\Ls_i}}$  depends on the type of admissible system. Let  $((\bb{p}_1,\bb{p}_2), (\bm{\phi}_1,\bm{\phi}_2))$ be an arbitrary $[\LL]$-admissible system and let $((\bb{p}_1,\bb{p}_2), (\bm{\iota}_1,\bm{\iota}_2))$ be the $\bb{K}$-fundamental system. In this proof, for clarity, we denote by $\bm{\Phi}_i$ the projection $\mathrm{proj}_{\mathcal{N}_{\Ls_i}}$ associated to $(\bb{p}_i,\bm{\phi}_i)$ and by  $\bb{I}_i$ the projection $\mathrm{proj}_{\mathcal{N}_{\Ls_i}}$ associated to $(\bb{p}_i,\bm{\iota}_i)$. Based on \eqref{eq:4.2.79}, we define the two norms  $\norm{\cdot}_{\M_{\LL}}$ and $\norm{\cdot}_{\tilde{\M}_{\LL}}$ as, $\forall f\in\M_{\LL}(\bb{K})$,
\begin{align}
    \norm{f}_{\M_{\Ls_1\otimes\Ls_2}}:=&\norm{[\Ls_1\otimes\Ls_2]\{f\}}_{\M(\R^2)}+\norm{[\bb{I}_1\otimes \Ls_2]\{f\}}_{1,\M(\R)}\nonumber\\
    &+\norm{[\Ls_1\otimes\bb{I}_2]\{f\}}_{\M(\R),1}+\norm{[ \bb{I}_1\otimes\bb{I}_2]\{f\}}_{\Yy^{\star}}
\end{align}
and
\begin{align}  \norm{f}_{\tilde{\M}_{\Ls_1\otimes\Ls_2}}:=&\norm{[\Ls_1\otimes\Ls_2]\{f\}}_{\M(\R^2)}+\norm{[\bm{\Phi}_1\otimes \Ls_2]\{f\}}_{1,\M(\R)}\nonumber\\
&+\norm{[\Ls_1\otimes\bm{\Phi}_2]\{f\}}_{\M(\R),1}+\norm{[ \bm{\Phi}_1\otimes\bm{\Phi}_2]\{f\}}_{\Yy^{\star}}.
\end{align}
Let $f$ be as in \eqref{eq:4.5.51}. Then, we define $\tilde{f}$ as 
\begin{align}
    \tilde{f}(\cdot)=&\int_{\bb{K}}g_{\bm{\phi}_1}(\cdot,x_1)g_{\bm{\phi}_2}(\cdot,x_2)\mathrm{d}m(x_1,x_2)+\sum_{n=1}^{N_1}p_{1,n}(\cdot)\left(\int_{\mathrm{K}_2}g_{\bm{\phi}_2}(\cdot,x_2)\mathrm{d}m_{1,n}(x_2)\right)\nonumber\\
    &+\sum_{n=1}^{N_2}\left(\int_{\mathrm{K}_1}g_{\bm{\phi}_1}(\cdot,x_1)\mathrm{d}m_{2,n}(x_1)\right) p_{2,n}(\cdot)+\sum_{n,n'=1,1}^{N_1,N_2}d_{n,n'}p_{1,n}(\cdot)p_{2,n'}(\cdot)\label{eq:A.7.113}
\end{align}
and observe that 
\begin{equation}
    \norm{f}_{\M_{\LL}}=\norm{\tilde{f}}_{\tilde{\M}_{\LL}}\leq\norm{f}_{\tilde{\M}_{\LL}}+\norm{f-\tilde{f}}_{\tilde{\M}_{\LL}}.
\end{equation}
Next, based on the calculations in Proposition \ref{prop:CompactRepresentation} we find that
\begin{align}
    f(t_1,t_2)-\tilde{f}(t_1,t_2)=&\sum_{n'=1}^{N_2}p_{2,n'}(t_2)\Ls^{-1}_1\left\{m_{2,n'}^{\mathrm{A}}\right\}(t_1)+\sum_{n=1}^{N_1}p_{1,n}(t_1)\Ls^{-1}_2\left\{m_{1,n}^{\mathrm{A}}\right\}(t_2)\nonumber\\
    &+\sum_{n,n'=1,1}^{N_1,N_2}(-d^{\mathrm{A}}_{n,n'}+d^{\mathrm{B}}_{n,n'}+d^{\mathrm{C}}_{n,n'})p_{1,n}(t_1)p_{2,n'}(t_2)\nonumber\\
    =&\underbrace{\sum_{n'=1}^{N_2}p_{2,n'}(t_2)\Ls^{-1}_{\bm{\phi}_1}\left\{m_{2,n'}^{\mathrm{A}}\right\}(t_1)}_{\Delta_{\RNum{2}}}+\underbrace{\sum_{n=1}^{N_1}p_{1,n}(t_1)\Ls^{-1}_{\bm{\phi}_2}\left\{m_{1,n}^{\mathrm{A}}\right\}(t_2)}_{\Delta_{\RNum{3}}}\nonumber\\
    &+\underbrace{\sum_{n,n'=1,1}^{N_1,N_2}(d^{\mathrm{A}}_{n,n'}+d^{\mathrm{B}}_{n,n'}+d^{\mathrm{C}}_{n,n'})p_{1,n}(t_1)p_{2,n'}(t_2)}_{\Delta_{\RNum{4}}}.
\end{align}
We further calculate that, for 
\begin{align}
    \forall i\in\{1,2\},\quad C_{i}&=N_i\underset{n\in[1\ldots N_i],x_i\in\mathrm{K}_i}{\text{sup}}\vert (g_{\Ls_i}^{\star}\ast\phi_{i,n})(x_i)\vert,
\end{align}
one has the upper-bounds

\begin{align}
    \norm{\Delta_{\RNum{2}}}_{\tilde{M}_{\LL}}&=\sum_{n'=1}^{N_2}\norm{m_{2,n'}^{\mathrm{A}}}_{\M(\R)}\leq C_{1}\norm{m}_{\M(\R^2)},\nonumber\\
    \norm{\Delta_{\RNum{3}}}_{\tilde{M}_{\LL}}&=\sum_{n=1}^{N_1}\norm{m_{1,n}^{\mathrm{A}}}_{\M(\R)}\leq C_{2}\norm{m}_{\M(\R^2)},\nonumber\\
\norm{\Delta_{\RNum{4}}}_{\tilde{M}_{\LL}}&\leq C\norm{\Delta_{\RNum{4}}}_1\leq C\left( C_{1}C_{2}\norm{m}_{\M(\R^2)}+\sum_{i=1}^{2}C_i\sum_{n=1}^{N_i}\norm{m_{i,n}}_{\M(\R)}\right),
\end{align}

where we used the fact that all norms are equivalent on $\R^{N_1}\otimes\R^{N_2}$, with $\norm{\cdot}_{\Yy^{\star}}\leq C\norm{\cdot}_1$. Consequently,
\begin{equation}
    \norm{f}_{\M_{\LL}}\leq\left(1+CC_{1}C_{2}+(1+C)C_{1}+(1+C)C_{2}\right)\norm{f}_{\tilde{\M}_{\LL}}.
\end{equation}
The proof of the existence of a constant $\tilde{C}$ such that 
$\norm{f}_{\tilde{\M}_{\LL}}\leq\tilde{C}\norm{f}_{\M_{\LL}}$ is similar and thus omitted. The more general proof where $((\tilde{\bb{p}}_1,\tilde{\bb{p}}_2), (\bm{\phi}_1,\bm{\phi}_2))$ is an admissible system with $(\tilde{\bb{p}}_1,\tilde{\bb{p}}_2)$ different from $(\bb{p}_1,\bb{p}_2)$ is the same, up to an additional change of basis. 
\end{proof}

\subsection{Proofs of Section 6.2 and 6.3.}
\label{app:5.3}

\begin{proof}[\textbf{Proof of Theorem \ref{th:RTlocalized}.}]
\hfill\\
\textbf{Step 1.}
It follows from Proposition \ref{prop:localizedjustification} that
\begin{equation}
\label{app:A.8.120}
    \forall x_i\notin[\text{min}(\phi_i^-,K_i^-),\text{max}(\phi_i^+,K_i^+)],\forall t_i\in[K_i^-,K_i^+]:\quad g_{\bm{\phi}_i}(t_i,x_i)=0.
\end{equation}
Finally, assume by contradiction that an extreme-point minimizer $f$ of the form \eqref{eq:5.2.66} has a knot $(x_{1,k},x_{2,k})\notin\tilde{\bb{K}}$ with associated weight $a_k\neq 0$. Then, 
\begin{align}
    \langle a_kg_{\bm{\phi}_1}(\cdot,x_{1,k})g_{\bm{\phi}_2}(\cdot,x_{2,k}),\bm{\nu}\rangle=\bb{0},\label{app:A.8.121}\quad\text{and}\\
    \snorm{a_kg_{\bm{\phi}_1}(\cdot,x_{1,k})g_{\bm{\phi}_2}(\cdot,x_{2,k})}_{\M_{\LL}}=\vert a_k\vert>0.\label{app:A.8.122}
\end{align}
The combination of \eqref{app:A.8.121} and \eqref{app:A.8.122} yields the new function 
\begin{equation}
   \tilde{f}(\cdot)=f(\cdot)-a_kg_{\bm{\phi}_1}(\cdot,x_{1,k})g_{\bm{\phi}_2}(\cdot,x_{2,k}),\quad\text{s.t.}\quad\mathcal{J}(\tilde{f})<\mathcal{J}(f).
\end{equation}
This is a contradiction. In turn, the inclusions  $z_{n',m'}\in\mathrm{K}_1$ and $y_{n,m}\in\mathrm{K}_2$ are proved likewise.
\hfill\\
\textbf{Step 2.} It follows from Step 1 that, if the system is $\bb{K}$-localized, then the extreme points of 
$\mathcal{V}$ are in $\M_{\LL}(\bb{K})$. Since $\M_{\LL}(\bb{K})$ is $\w$-closed, the $\w$ closure of these extreme points (i.e., $\mathcal{V}$) must be a subset of $\M_{\LL}(\bb{K})$. Therefore, $\mathcal{V}\subset\M_{\LL}(\bb{K}).$
\end{proof}

\begin{proof}[\textbf{Proof of Corollary \ref{coro:interior}.}]
Consider a generic extreme-point solution $f^{\star}$ with knots $(x_{1,k},x_{2,k})$, $ y_{n,m},$ and $z_{n',m'}$. 
\begin{itemize}
    \item If there exists a knot $(x_{1,k},x_{2,k})$, with associated weight $a_k\neq0$, that is such that $x_{1,k}=K_{1}^+$ or $x_{2,k}=K_{2}^+$, then $\forall\bb{t}=(t_1,t_2)\in\bb{K},$
    \begin{align}
       & a_kg_{\Ls_1}(t_1-x_{1,k})g_{\Ls_2}(t_2-x_{2,k})=0\nonumber\\
       \Rightarrow\hspace{0.1cm}&\langle a_kg_{\Ls_1}(t_1-x_{1,k})g_{\Ls_2}(t_2-x_{2,k}),\bm{\nu}\rangle=\bb{0}.
    \end{align}
    This yields a contradiction with the minimality of $f^{\star}$, since the new function
    \begin{equation}
        f^{\star\star}(\bb{t})=f^{\star}(\bb{t})-a_kg_{\Ls_1}(t_1-x_{1,k})g_{\Ls_2}(t_2-x_{2,k})
    \end{equation} has the same measurement value but a lower regularization.
    
    \item If there exists a knot $(x_{1,k},x_{2,k})$, with associated weight $a_k\neq0$, that is such that $x_{1,k}=K_1^-$ and $x_{2,k}\in]K_{2}^-,K_2^+[$, then $\forall\bb{t}\in\bb{K},$
    \begin{align}
        a_kg_{\Ls_1}(t_1-x_{1,k})g_{\Ls_2}(t_2-x_{2,k})=&a_k\frac{(t_1-K_1^-)_+^{N_1-1}\e^{\alpha_1(t_1-K_1^-)}}{(N_1-1)!}g_{\Ls_2}(t_2-x_{2,k})\nonumber\\
        =&a_k\frac{(t_1-K_1^-)^{N_1-1}\e^{\alpha_1(t_1-K_1^-)}}{(N_1-1)!}g_{\Ls_2}(t_2-x_{2,k})\nonumber\\
        =&a_kp_{1,N_1}(t_1)g_{\Ls_2}(t_2-x_{2,k}).
    \end{align}
    Therefore, the new function
    \begin{equation}
        f^{\star\star}(\bb{t})=f^{\star}(\bb{t})-a_kg_{\Ls_1}(t_1-K_1^-)g_{\Ls_2}(t_2-x_{2,k})+a_kp_{1,N_1}(t_1)g_{\Ls_2}(t_2-x_{2,k})
    \end{equation}
    is such that
    \begin{equation}
        \forall\bb{t}\in\bb{K}, f^{\star}(\bb{t})= f^{\star\star}(\bb{t}),\quad\text{and}\quad\snorm{f^{\star}}_{\M_{\LL}}=\snorm{f^{\star\star}}_{\M_{\LL}}.
    \end{equation}
    \item The argumentation for when $x_{2,k}=K_2^-$ and $x_{1,k}\in]K_{1}^-,K_1^+[$ or, for when $x_{1,k}=K_1^-$ and $x_{2,k}=K_2^-$, is similar and omitted.
\end{itemize}
On the 3 previous bullet points, the iteration on the knots $(x_{1,k}, x_{2,k})$ that do not belong to the interior $\overset{\circ}{\bb{K}}$ yields a new extreme-point solution $f^{\star\star}$ that satisfies Items 2 and 3, and such that all of its knots $(x_{1,k}, x_{2,k})$ are in $\overset{\circ}{\bb{K}}$. Finally, one can find an $f^{\star\star}$ that satisfies Items 2, 3 with knots $y_{n,m}\in\overset{\circ}{\mathrm{K}}_1$ and $z_{n',m'}\in\overset{\circ}{\mathrm{K}}_2$, by a similar argumentation.
    
\end{proof}

\bibliographystyle{elsarticle-harv}
\bibliography{ref}

\end{document}